\documentclass[%
 aip,
 amsmath,amssymb,
 reprint,%
]{revtex4-1}

\usepackage{graphicx}% Include figure files
\usepackage{dcolumn}% Align table columns on decimal point
\usepackage{bm}% bold math
\usepackage[utf8]{inputenc}
\usepackage[T1]{fontenc}
\usepackage{mathptmx}
\usepackage{etoolbox}
\usepackage{color,soul}
\usepackage{svg}
\usepackage{listings}
\usepackage{changes}
\usepackage{float}
\usepackage{todonotes}
\usepackage[hidelinks]{hyperref}
\usepackage{framed}
\usepackage{braket}
\usepackage{changepage}

\usepackage{amsfonts}

\usepackage[dvipsnames]{xcolor}
\usepackage{tikz}
\usetikzlibrary{positioning,decorations.pathreplacing,arrows.meta,shapes.geometric,calc}

\definecolor{codebg}{RGB}{248,248,248}
\definecolor{codecomment}{RGB}{80,130,80}
\definecolor{codekeyword}{RGB}{35,70,140}
\definecolor{codestring}{RGB}{160,60,60}
\definecolor{codenumber}{RGB}{120,120,120}

\lstdefinestyle{pythonstyle}{
  language=Python,
  basicstyle=\ttfamily\small,
  backgroundcolor=\color{codebg},
  keywordstyle=\color{codekeyword}\bfseries,
  commentstyle=\color{codecomment}\itshape,
  stringstyle=\color{codestring},
  numbers=left,
  numberstyle=\tiny\color{codenumber},
  numbersep=8pt,
  frame=single,
  rulecolor=\color{codenumber},
  framesep=5pt,
  breaklines=true,
  breakatwhitespace=true,
  showstringspaces=false,
  keepspaces=true,
  columns=fullflexible,
  tabsize=4,
  captionpos=t,
  float=t,
  aboveskip=0.8em,
  belowskip=0.8em,
  xleftmargin=1.5em,
  framexleftmargin=1.2em
}

\definecolor{lightergray}{gray}{0.95} %cadre
\definecolor{orangebrique}{RGB}{230,100,4} %def
\definecolor{emeraude}{RGB}{0,150,120} %prop
\definecolor{joliejaune}{RGB}{220,140,15}
\tikzstyle{block} = [draw, rectangle, 
minimum height=1em, minimum width=6em]

\makeatletter
\def\@email#1#2{%
 \endgroup
 \patchcmd{\titleblock@produce}
  {\frontmatter@RRAPformat}
  {\frontmatter@RRAPformat{\produce@RRAP{*#1\href{mailto:#2}{#2}}}\frontmatter@RRAPformat}
  {}{}
}%
\makeatother

\def\ie{i.e.}
\def\eg{e.g.}

\definecolor{formalshade}{rgb}{0.95,0.95,1}
\definecolor{darkblue}{rgb}{0.0, 0.0, 0.55}
\definecolor{tabblue}{RGB}{31, 119, 180}
\definecolor{tabred}{RGB}{214, 39, 40}

\newenvironment{formal}{%
  \MakeFramed{\advance\hsize-\width\FrameRestore}%
  \noindent\hspace{-4.55pt}% disable indenting first paragraph
  \begin{adjustwidth}{}{7pt}%
  \vspace{2pt}\vspace{2pt}%
}
{%
  \vspace{2pt}\end{adjustwidth}\endMakeFramed%
}

\begin{document}

%\preprint{AIP/123-QED}

\title[Analog and hybrid computing in the 21st century]{Analog and hybrid computing in the 21st century}
\author{Bernd Ulmann}
\thanks{Corresponding authors: ulmann@anabrid.com, lminati@buaa.edu.cn, and l.boettcher@fs.de}
\affiliation{FOM University of Applied Sciences for Economics and Management, 60323 Frankfurt a.\ M., Germany}
\affiliation{anabrid GmbH, 60486 Frankfurt a.\ M., Germany}

\author{Olivier Bournez}
\affiliation{LIX, CNRS, Ecole Polytechnique, Institut Polytechnique de Paris, 91120 Palaiseau, France}

\author{Corey J. Maley}
\affiliation{Department of
Philosophy, Purdue University, West
Lafayette, IN 47907, United States of America}

\author{Zhong Sun}
\affiliation{Institute for Artificial Intelligence, School of Integrated Circuits and Beijing Advanced Innovation Center for Integrated Circuits, Peking University, Beijing 100871, China}

\author{Grace Li Zhang}
\affiliation{Hardware for Artificial Intelligence Group, Technical University of Darmstadt,
64283 Darmstadt, Germany}

\author{Paul Manea}
\affiliation{Peter Grünberg Institut (PGI-14),
Forschungszentrum Jülich GmbH, 52068 Aachen, Germany}

\author{Michael D\"{o}ll}
\affiliation{Department of Computer Science, Darmstadt University of Applied Sciences, 64295 Darmstadt, Germany}

\author{Frank Delvigne}
\affiliation{Gembloux Agro-Bio Tech, Université de Liège, Gembloux, 5030, Belgium}

\author{Alex J.H. Fedorec}
\affiliation{Division of Biosciences, University College London, London, WC1E 6BT, United Kingdom}

\author{Daniel Garza}
\affiliation{Universit\'{e} Paris-Saclay, PROSE, INRAE, Antony, 92160, France}

\author{Lucas Wetzel}
\affiliation{anabrid GmbH, 60486 Frankfurt a.\ M., Germany}

\author{Sara Obergassel}
\affiliation{anabrid GmbH, 60486 Frankfurt a.\ M., Germany}

\author{Pavol Pitonak}
\affiliation{anabrid GmbH, 60486 Frankfurt a.\ M., Germany}

\author{Christopher Jung}
\affiliation{anabrid GmbH, 60486 Frankfurt a.\ M., Germany}

\author{Daniel Thuerck}
\affiliation{anabrid GmbH, 60486 Frankfurt a.\ M., Germany}

\author{Mason A. Porter}
\affiliation{Nordita (Nordic Institute for Theoretical Physics), KTH Royal Institute of Technology, Stockholm University, and Uppsala University, Sweden; Department of Mathematics, Stockholm University, Sweden;
Department of Mathematics, University of California, Los Angeles; United States of America; Department of Sociology, University of California, Los Angeles; United States of America; and Santa Fe Institute, United States of America}

\author{Marcelo J. Rozenberg}
\affiliation{
Universit\'e Paris-Cit\'{e}, CNRS, Integrative Neuroscience and Cognition Center, 75006, Paris, France}
\affiliation{Department of Physics, University of California, San Diego, CA 92093, USA}

\author{Frederic Bergner}
\affiliation{German Aerospace Center (DLR), Institute of Aerodynamics and Flow Technology, 37073 G\"{o}ttingen, Germany}
\affiliation{University of G\"{o}ttingen, Institute for Numerical and Applied Mathematics, 37073 G\"{o}ttingen, Germany}

\author{Bhaskar Choubey}
\affiliation{Institute of Analogue Circuits and Image Sensors, Universit\"{a}t Siegen, 57076 Siegen, Germany}

\author{Manyu Zhao}
\affiliation{Faculty of Engineering, Free University of Bozen-Bolzano, 39100 Bozen-Bolzano, Italy}

\author{Ludovico Minati}
\thanks{Corresponding authors: ulmann@anabrid.com, lminati@buaa.edu.cn, and l.boettcher@fs.de}
\affiliation{International Research Center of Complexity Sciences, Hangzhou International Innovation Institute, Beihang University, 311115 Hangzhou China}
\affiliation{School of Life Science and Technology, University of Electronic Science and Technology of China, 611731 Chengdu, China}

\author{Lucas B\"{o}ttcher}
\thanks{Corresponding authors: ulmann@anabrid.com, lminati@buaa.edu.cn, and l.boettcher@fs.de}
\affiliation{Department of Computational Science and Philosophy, Frankfurt School of Finance \& Management, 60322 Frankfurt a.\ M., Germany}

\date{\today}% It is always \today, today,

\begin{abstract}
    Analog computers solve problems through physical representations of mathematical variables and relations, often using the dynamics of electronic circuits and their associated signals. Many analog computers directly represent nonlinear dynamical systems; they encode their solutions in the signals that they generate.
    Although most modern scientific computing relies on digital architectures with inherently discrete-valued representations and algorithmic approaches, the design constraints of such
    architectures increasingly limit both performance and energy efficiency. Analog systems can circumvent some of these constraints; they execute specific computations very efficiently by mapping them directly onto physical processes using smaller and simpler computing elements than the ones in digital architectures. 
    In this article, we provide an interdisciplinary perspective on modern analog and hybrid (\ie, analog-digital) computing. We discuss how modern scientific computing workflows can leverage analog approaches in various applications, such as in simulating dynamical systems, performing optimization, and enabling more energy-efficient implementations of artificial intelligence.
\end{abstract}

\maketitle

%%%%

\begin{quotation}

    Given the ever-increasing demand for energy-efficient and scalable computing, there has been a resurgence of interest in analog computers. Historically, mechanical analog computers were used to compute celestial positions, predict tides, and solve differential equations.\cite{ulmann2026analog}
    Such systems are not suitable for modern scientific-computing workflows, but the underlying idea of configuring a computing system to the structure of a problem remains central to many contemporary scientific and industrial domains. 
    Modern analog computers follow this paradigm by solving problems using physical processes and magnitudes for both computation and representation. This contrasts with general-purpose digital computing, which represents information using discrete values and commonly relies on von Neumann architectures,\cite{patterson2016computer} where data movement between physically separated processing and memory increasingly limits performance and energy efficiency. In this article, we provide an interdisciplinary perspective on state-of-the-art analog and hybrid computing, highlight the connections between analog computation and nonlinear dynamics, and discuss how analog approaches can benefit the simulation of dynamical systems, the solution of optimization problems, and the implementation of more energy-efficient artificial-intelligence systems. We also consider how to integrate analog computers into heterogeneous computing stacks.
\end{quotation}
\section{\label{sec:intro}Introduction}
What is an analog computer, and how does it differ from our now-ubiquitous digital computers? A digital computer, which is more aptly called a ``stored-program digital computer'', is a machine with a fixed internal structure that solves problems by executing sequences of simple instructions that operate on data that is represented digitally using binary values.\cite{patterson2016computer} This concept is incredibly powerful, as such a machine can carry out any computation that is describable by an algorithm if there are sufficient time and memory resources.
It rests on both Turing's universal-machine construction and the Church--Turing thesis, which identifies algorithmic solvability with computability by a Turing machine.\cite{Tur36,Church36} We will return to this thesis and to its analog counterpart in Sec.~\ref{sec:new_quest}. 

The problem with current digital computers is that they are about to hit fundamental physical boundaries in integration density (\ie, the number of transistors per unit area of an integrated circuit), execution speed, and power consumption.\cite{Horowitz2014} 

In the context of integration density, the term ``technology node'' historically referred approximately to a characteristic feature size that is measured in nanometers (nm). However, semiconductor manufacturers now use node names such as ``5 nm'' and ``2 nm'', even though these numbers no longer correspond directly to any specific physical dimensions of underlying integrated circuits.\cite{badaroglu2021}
As critical device and interconnect dimensions approach the nanometer regime, heat dissipation and quantum effects impose increasingly severe constraints. Accordingly, continued increases in integration density will become both more difficult and more costly.

Clock frequencies of digital computers have essentially stalled in the few-gigahertz (GHz) range over the last two decades, primarily  because the energy of digital circuits increases faster than linearly with clock frequencies.\cite{zyuban2004integrated} 
In practice, doubling the clock rate from 3~GHz to 6~GHz typically results in more than twice the power consumption. This is not only inefficient but also a technological challenge, as the electrical power that is consumed by a chip is converted into heat, which then must be dissipated from a rather small chip area.\cite{esmaeilzadeh2012dark,Horowitz2014}

The current demand for computing power, which is driven in part by the rapid growth of artificial intelligence (AI), is leading to increasingly large data centers with substantial energy consumption. This has sparked discussions about recommissioning old nuclear power plants or building new nuclear reactors that are dedicated to individual data centers.\cite{hemphill2024us} Although these are interesting from a technological point of view, their widespread deployment entails substantial economic and environmental challenges.\cite{ruble2026us}
Accordingly, to satisfy current and future demand for computing power, it is necessary to explore different computational paradigms beyond the usual digital computing. One of these paradigms is ``analog computing'', which is the focus of our perspective article.

An analog computer differs fundamentally from a digital computer. First, in an analog computer, values are typically not represented by binary digits; instead, they are encoded by magnitudes of physical processes (\eg, by voltages or currents). 
Second, an analog computer consists of a potentially large number of computing elements, which each implement a specific basic mathematical operation (such as summation, multiplication, and time integration). One then connects these computing elements with each other, thereby implementing an electronic model (\ie, an ``analogue'') of a problem that one seeks to solve. This model-based approach to computing gave rise to the name ``analog computing".\footnote{Unfortunately, many people perceive the notion of ``analog computing'' as being old-school, slow, or otherwise suboptimal. Accordingly, it may be worthwhile to abandon the name ``analog computer'' for future discussions and replace it by terms such as ``natural computing'' or ``physical computing''.}
\begin{figure}
    \begin{tikzpicture}[
    line/.style={draw,line width=0.8pt},
    arrow/.style={-{Latex[length=3mm,width=2mm]},line width=0.9pt},
    circ/.style={draw,circle,line width=0.8pt,minimum size=1.0cm,inner sep=0pt},
    every node/.style={font=\normalsize},
    branchdot/.style={circle,fill,inner sep=1.55pt}
]

% Main block: thin rectangle plus triangular wedge
\draw[line] (0,0) rectangle (0.50,2.00);
\draw[line] (0.50,0.00) -- (2.00,1.00) -- (0.50,2.00);

% Initial condition circle above the block
\node[circ] (init) at (0.25,3.00) {$x(0)$};
\draw[line] (0.25,2.00) -- (init.south);
\draw[line] (init.north) -- (0.25,3.95);
\node[above=-1pt] at (0.25,3.95) {$-1$};

% Output of the triangular block
\coordinate (out) at (2.50,1.00);
\draw[line] (2.00,1.00) -- (out);
\node at (2.06,1.43) {$x(t)$};

% Branch going upward to displayed output x
\draw[line] (out) -- (2.50,2.50);
\draw[arrow] (2.50,2.50) -- (3.25,2.50);
\node[right=3pt] at (3.25,2.50) {$x(t)$};

% Multiplication by lambda
\node[circ] (lambda) at (3.50,1.00) {$\lambda$};
\draw[line] (out) -- (lambda.west);

% Feedback path carrying -xdot
\draw[line]
    (lambda.east) -- (4.45,1.00) -- (4.45,-0.50)
    -- (-0.50,-0.50) -- (-0.50,1.00) -- (0,1.00);

\node[left=0.12cm] at (0,1.36) {$-\dot{x}$};

% Solid branch dot, drawn last
\node[branchdot] at (out) {};

\end{tikzpicture}
    \caption{An analog computer program for solving the ordinary differential equation (ODE) $\dot{x}(t) = -\lambda x(t)$. The symbol with a triangle and rectangle together denotes an integrator, which performs a continuous-time integration over the sum of its inputs. (In this example, we use only one input.) For technical reasons, integrators perform an implicit sign conversion. The circle on top of the integrator gives its initial condition $x(0) = 1$. For this initial condition and times $t\geq 0$, the integrator output represents the state $x(t)$. In the feedback loop, $x(t)$ is multiplied by $\lambda$ and applied to the integrator input.}
    \label{pic_decay_example}
\end{figure}
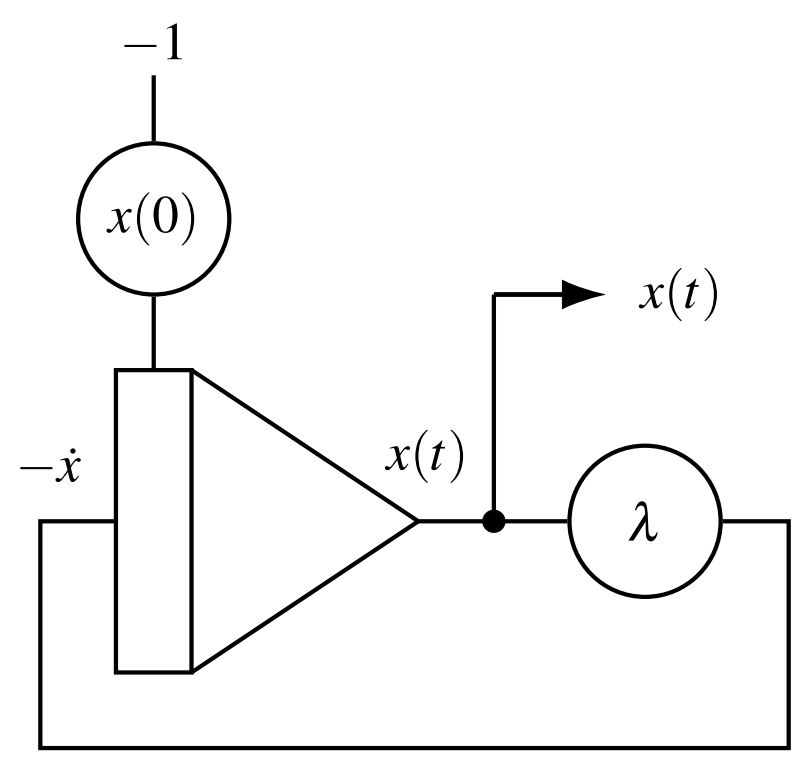
\begin{lstlisting}[
  style=pythonstyle,
  float=t,
  captionpos=t,
  caption={Digital computer program for solving the ODE $\dot{x}(t) = -\lambda x(t)$ with the initial condition $x(0) = 1$ using the explicit Euler method.},
  label={lst:euler}
]
lambda_ = 1.0
dt = 0.01
x = 1.0

for n in range(num_steps):
    x_dot = -lambda_ * x
    x = x + dt * x_dot
\end{lstlisting}

What is it that actually makes analog computing a notion of \textit{computing}? What makes it similar enough to the digital-computing paradigm to count as \textit{computation}, yet different enough to make it a distinct type of computation? In other words, what is computation in general, such that digital computing is one type and analog computing is another type?\footnote{We are going to entirely avoid discussing quantum computing in this context.} The answer that seems to best fit this perspective is that computation involves the following items:
\begin{itemize}
\itemsep-0.25em
    \item{physical representations or variables (which can represent physical systems or mathematical systems), which}
    \item{are manipulated by physical mechanisms}
    \item{and change the values of those representations or variables by manipulating the physical properties that are responsible for representing them.}
\end{itemize}
Digital computers use digital representations, which are manipulated by electronic mechanisms in terms of discrete voltages that are controlled by a sequentially executed algorithm. Analog computers do not employ such algorithms; 
instead, they rely on setting up a model. Electronic analog computers are mechanisms that manipulate continuous representations in terms of voltages.\nocite{maley_analogue_2023}
\footnote{See Ref.~\onlinecite{maley_analogue_2023} for a philosophical account of analog computing.}

An analog computer program is a directed, cyclic graph that specifies the connections between computing elements. In Fig.~\ref{pic_decay_example}, we show a program for solving the ordinary differential equation 
(ODE) $\dot{x}(t) = -\lambda x(t)$ with the initial condition $x(0) = 1$. In contrast with this analog approach, in Listing~\ref{lst:euler}, we show a digital solution of the same initial-value problem using the explicit Euler method.

Integrators are the only computing elements with time-dependent behavior. An integrator has three operation modes: IC (initial condition), OP (operate), and HALT. In the IC mode, one sets the output of the integrator to its initial condition. In the OP mode, the integrator performs an integration with respect to time. One can temporarily stop the integration using the
HALT mode, which is useful to read out values. The operation speed of analog-computer integrators determine the 
time that it takes to obtain a solution. Let $k_0$ denote the time-scaling factor of an integrator.
The value $k_0 = 1$ signifies that integrating the value $-1$ for one second with the initial value $0$ yields an output value of $1$. 
An integrator with a time-scaling factor $k_0 = 10$ runs the same integration in $0.1$ physical seconds. For modern analog computers, $k_0$ ranges from approximately $10^2$ to approximately $10^4$ for anabrid GmbH's LUCIDAC and REDAC systems,\footnote{The name LUCIDAC combines ``lucid'' (clear, software-programmable patching) with ``DAC'' (digital-to-analog converter). The acronym REDAC stands for ``reconfigurable discrete analog computer''. The LUCIDAC and REDAC systems are both fully reconfigurable modern hybrid computer systems. The LUCIDAC system is the smaller of the two.} whereas recent {analog computing systems that are implemented on a single integrated circuit} can exceed $k_0 = 10^6$. 
With such large time-scaling factors, it takes in the range of microseconds for a single integration of coupled ODEs and spatially discretized PDEs.\footnote{For boundary-value problems, one needs to run a circuit many times.}
The high speed of analog computers in these tasks results primarily from the fact that there is no memory from which to read instructions to execute sequentially. Instead, all computing elements of a program work in a fully parallel way.

In addition to their fast solution times, analog computers can achieve very high energy efficiency. Representing values as voltages or currents can reduce the energy overhead that is associated with digital switching operations and data movement.

The aim of the present perspective article is to provide an interdisciplinary overview of modern analog and hybrid computing. In this article, we bring perspectives across a wide range of fields, including physics, mathematics, computer science, electrical engineering, biology, and philosophy. The individual contributions are guided by the following overarching questions:
\begin{itemize}
\itemsep-0.25em 
\item{What can analog computers compute, and how does analog computation compare with digital computation?}
\item{What are the current paradigms of analog computing, and which applications can they address?}
\item{How can one design general-purpose analog computers (GPACs) and make them accessible to both academia and industry?}
\item{What are relevant scientific and industrial application domains, and can analog computers help address the energy challenges that are associated with modern artificial intelligence?}
\end{itemize}

We will address these questions and related questions in the four main parts of our perspective. In Sec.~\ref{sec:new_quest}, we discuss the new quest for analog computing by (i) reviewing results about computability and complexity and by (ii) establishing connections between GPACs, polynomial ODEs, and other models of computation. In Sec.~\ref{sec:current_paradigms}, we overview current analog-computing paradigms, which range from in-memory computing and biological computing to modern programmable GPAC hardware. 
Building on this overview of hardware platforms, in Sec.~\ref{sec:software}, we discuss the software stacks and system-development tools that one needs to ensure that programmable GPAC systems are accessible for scientific and industrial applications. In Sec.~\ref{sec:applications}, we overview several application areas that are suitable for analog computers. These application areas include dynamical-systems modeling, PDEs, optimization, edge computing, computing based on minimalist chaotic circuits, and artificial intelligence.
 In Sec.~\ref{sec:summary}, we summarize our perspective's key topics. We then conclude in Sec.~\ref{sec:conclusion} by posing several open research questions.\footnote{In our present perspective, we do not cover the rich and interesting history of analog computing. For a comprehensive account, see Ref.~\onlinecite{ulmann2026analog}.}
\section{\label{sec:new_quest}The new quest for analog computing---By
Olivier Bournez}
\subsection{A mathematical model of analog computers}
The analog computers described in Sec.~\ref{sec:intro} raise basic theoretical questions: what exactly can such machines compute and how does this compare to a digital computer? A first answer to these questions goes back to Shannon, who in 1941 introduced a mathematical idealization of Bush's differential analyzer, the general purpose analog computer (GPAC).\cite{Sha41,Bus31} A GPAC is built from the same primitives as a physical analog computer: constant units, adders, multipliers, and integrators wired together (see Fig.~\ref{pic_decay_example}). The natural object it produces is a function of time. One can show that the functions a GPAC generates are exactly the components of the solutions of polynomial ODEs of the form $y' = p(y)$, where $p$ is a vector of polynomials~\cite{GC03,InformationAndComputation2017}; making this statement precise took several decades, through successive corrections of Shannon's original claim.\cite{PER74,LR87} In other words, programming an analog computer in this sense amounts to choosing a polynomial vector field together with an initial condition, and the quantities it produces are generated in real time. The value at time $t$ is available at physical time $t$. Such quantities are in particular
differentially algebraic. They include polynomials, the exponential, sine and cosine, and functions built from them by addition, multiplication, and integration.
\subsection{Generating versus computing}
Shannon's original notion, which we call GPAC-generability, is more demanding than it first appears. There is no input other than time itself. The machine is started, and the observed signal must be $f(t)$ at physical time $t$, in real time. Two consequences follow. First, some perfectly ordinary functions fall outside its scope; the Gamma function $\Gamma(x)=\int_0^{\infty} u^{x-1}e^{-u}\,\mathrm{d}u$, for instance, is not GPAC-generable, since by a classical theorem of H\"older it is not
differentially algebraic.\cite{Hol87}
Second, real-time generation is not what we usually mean today by ``computing a function''. In modern terms, computing $f$ on input $x$ means producing, after some computation time $t(x)$, an approximation of $f(x)$ to any prescribed accuracy, with longer computation giving better precision. Encoding the input as an initial condition and reading the result off one coordinate of the trajectory as it converges yields a second, more permissive notion, GPAC-computability. Under this notion the Gamma function becomes computable,\cite{Gra04} a first sign that
GPAC-computability, rather than generability, is the right yardstick when comparing analog and digital machines.
\subsection{Analog models of computation are Turing complete}
\label{sec:analog_models_turing_complete}
\begin{figure}
    \centering
    \begin{tikzpicture}[
    xscale=0.75, yscale=0.9, auto, node distance=10cm,
    more/.style = {block, dotted, thick, fill=gray!8, text=gray!45!black,
                   align=center, text width=5.8em, font=\footnotesize}
  ]
    \node[block, fill=emeraude!10, yshift=0.5cm]        (tm)       {Turing Machines};
    \node[block, below=0.5cm of tm, fill=joliejaune!10] (prg)      {Programs in C};
    \node[block, below=0.5cm of prg, fill=orangebrique!10] (circuits) {Boolean Circuits};
    \node[block, below=0.5cm of circuits, fill=DarkOrchid!10] (lambda) {Algebraic Models};
    \node[more, below=0.5cm of lambda]                  (moreL)    {other discrete models};
    \node[block, below=1.5cm of moreL, fill=Violet!10, very thick] (quantum) {Quantum Models};

    \node[block, right=3cm of tm, fill=MidnightBlue!10] (gpac) {GPAC};
    \node[block, below=1cm of gpac, fill=RoyalBlue!10]  (pODE) {Polynomial ODEs};
    \node[block, below=1cm of pODE, fill=Cerulean!10]   (bio)  {Biochemical Systems};
    \node[more, below=1cm of bio]                       (moreR) {other continuous models};

    \draw[<->] (tm)       -- (prg);
    \draw[<->] (prg)      -- (circuits);
    \draw[<->] (circuits) -- (lambda);
    \draw[<->] (gpac)     -- (pODE);
    \draw[<->] (pODE)     -- (bio);

    \draw[dotted, thick] (lambda) -- (moreL);
    \draw[dotted, thick] (bio)    -- (moreR);

    \draw[<->, dashed] (moreL) -- (quantum);

    \coordinate (Lb) at ($(circuits.east |- tm.north)+(0.30,0)$);
    \coordinate (Rb) at ($(bio.west |- gpac.north)-(0.45,0)$);
    \draw[decoration={brace, raise=5pt}, decorate] (Lb) -- (Lb |- moreL.south);
    \draw[decoration={brace, raise=5pt}, decorate] (Rb |- moreR.south) -- (Rb);

    \draw[<->, dashed] ($(Lb |- circuits)+(0.45,0)$)
        -- node[below=6pt] {?} ($(Rb |- circuits)-(0.45,0)$);
\end{tikzpicture}
    \caption{Models of computation and the extended Church--Turing thesis. The Church--Turing thesis states that all reasonable and sufficiently powerful models of computation are equivalent in terms of computability. The effective, or extended, Church--Turing thesis asks whether this equivalence also holds at the complexity level, namely whether such models can simulate one another with at most polynomial overhead. Quantum computation is a possible boundary of this principle. While it does not change what is computable, it may change what is efficiently computable. The analogous question for analog computation is whether it remains within this complexity-level equivalence, or whether it constitutes another boundary of the Church--Turing thesis. In the case of polynomial ODEs, this question can be made precise through the length of trajectories. Polynomial-time computation corresponds to solutions of polynomial ODEs whose trajectories have polynomial length.
}
    \label{fig:figTM}
\end{figure}
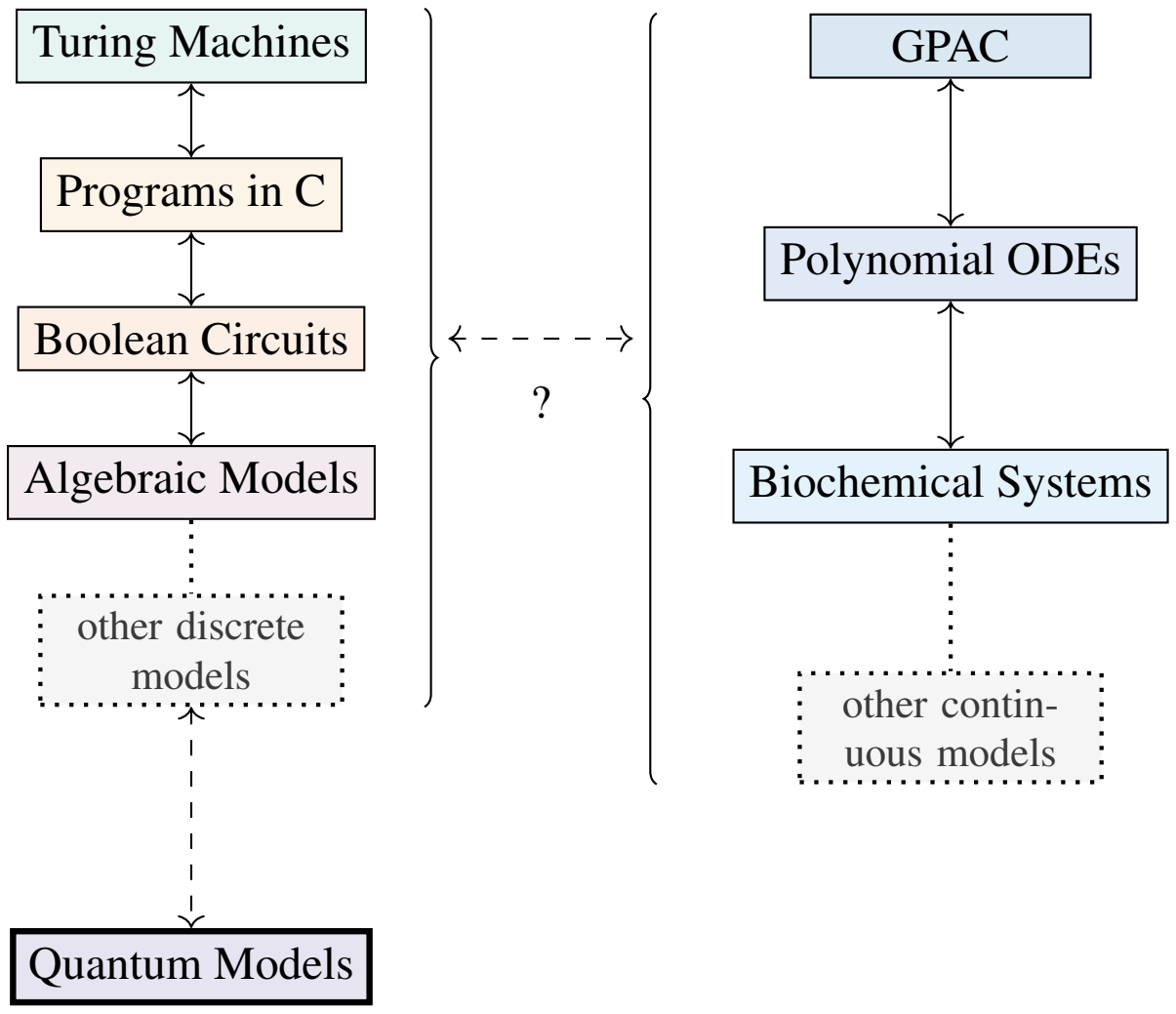
Analog and digital models have the same computational power, in the sense that anything one of them can compute, the other can compute too. In particular, a real function on a compact domain is computable by a Turing machine, in the sense of computable analysis, if and only if it is computable by a GPAC.\cite{JOC2007} Polynomial differential equations are thus neither weaker nor stronger than Turing machines.

This equivalence dissolves two opposite misconceptions. The first underestimates analog computers, treating them as nothing more than imprecise, noisy versions of digital ones, and therefore strictly weaker. However, in the converging sense, they compute exactly the Turing-computable functions. The second overestimates them. Because a single real number carries unbounded information, an idealized continuous device that reads and manipulates exact reals can, in principle, decide problems no Turing machine can. Recurrent neural networks with exact real weights, or wave equations with computable initial data whose solution is not computable, are the standard illustrations.\cite{SS94,PR81} This apparent hypercomputation is an artifact of unphysical idealizations, namely unbounded precision and exact, noise-free constants. It disappears as soon as the computation is required to be robust to arbitrarily small perturbations, the natural counterpart of finite precision, and what then remains is exactly Turing
power.\cite{BGHRobust10} The robust version also asks nothing exotic of the machine itself. A Turing machine can be simulated by a polynomial ODE with rational coefficients~\cite{JOC2007}, so no uncomputable real constant is ever needed.
\subsection{The effective Church--Turing thesis, and the question of speed}
The equivalence between analog and digital computation must be read carefully, because the Church--Turing thesis comes in several distinct forms. Its computability form, applied to physical devices and sometimes called the physical Church--Turing thesis, concerns what can be computed at all, and the result just stated is its analog counterpart. No robust analog machine computes anything beyond the reach of a Turing machine, and conversely. A stronger form, the effective or extended Church--Turing thesis, concerns what can be computed \emph{efficiently}, asserting that any reasonable physical model can be simulated by a Turing machine with only polynomial overhead. Asking whether analog machines can compute \emph{faster} is asking
whether they obey this second, stronger thesis.

The same question drives quantum computing, where the extended thesis is widely believed to fail, since quantum algorithms appear to offer super-polynomial speed-ups for problems such as integer factoring. The analog case is more subtle, and one point is decisive. A differential equation can always be sped up by rescaling time, exactly as raising the time-scaling factor $k_0$ of Sec.~\ref{sec:intro}
does on a real machine, so physical duration alone is not a meaningful measure of effort. A faithful measure is the length of the solution curve, the distance actually traveled along the trajectory, which is invariant under such rescaling. With this measure the equivalence extends to complexity: a function is computable in polynomial time if and only if it is computable by a polynomial differential equation whose solution has polynomial length.\cite{JournalACM2017} Additionally, it has recently been shown that the space, or memory, of a Turing machine corresponds to the precision at which the continuous quantities must be maintained.\cite{Icalp2024} Whether analog machines can genuinely compute faster then becomes a precise question about the length and precision of solutions of polynomial ODEs.

Beyond their conceptual interest, these unexpected links between ODEs and computation have produced new results in computer science in their own right, such as the construction of a universal ODE,\cite{LMCS2018} which answers a question of Rubel.\cite{Rub81} Refer to Ref.~\onlinecite{bournez2021survey} for a recent survey.
\subsection{Bridges to other natural models}
These results matter beyond electronics because the same mathematics governs other ways of computing with continuous quantities (see Fig.~\ref{fig:figTM}). The recurrent neural networks mentioned in Sec.~\ref{sec:analog_models_turing_complete} are one such family, with the caveat that their classical universality results concern \emph{discrete-time} networks, whereas differential analyzers and electronic integrator networks are genuinely continuous-time systems, for which the GPAC and polynomial ODEs are the appropriate model. The same polynomial vector fields also describe the kinetics of chemical reaction networks, and continuous chemical reaction networks are themselves Turing-complete, with GPAC-computability corresponding to computability by such
networks while preserving resources.\cite{CMSB17vantardise} One theory thus spans electronic analog computers, neural networks, and molecular or biochemical computation, including the biological systems discussed in
Sec.~\ref{sec:bio_computing}, which is what makes the GPAC a common foundation rather than a model tied to a single technology. A detailed, theoretical computer science-oriented can be found in Ref.~\onlinecite{bournez2021survey}.
\section{\label{sec:current_paradigms}Current paradigms of analog computing}
\subsection{\label{sec:memristive}In-memory computing---By Zhong Sun, Grace Li Zhang, Paul Manea, and Michael D\"oll}
Matrix computations are essential for a wide range of applications, including neural networks, optimization, scientific computing, and signal processing. In the big data era, however, the ever-growing matrix dimensions, combined with high computational complexity and the memory-wall bottleneck of conventional von Neumann architectures, pose substantial challenges for digital implementations. In-memory computing using analog devices such as memristors/resistive random-access memory (RRAM) has emerged as a promising solution that overcomes the energy and data movement limitations of conventional architectures, providing a substantial improvement in parallelism by leveraging the inherently parallel memory array structure.\cite{sun2023full} 
\subsubsection{Matrix--vector multiplication and matrix equation solving}
Resistive random-access memory crossbar arrays can efficiently perform in-memory matrix--vector multiplications (MVMs) for neural network inference and training, which is useful in edge artificial intelligence applications. Furthermore, by configuring closed-loop connections between the RRAM crossbar and peripheral analog circuits (\eg, operational amplifiers, transimpedance amplifiers), the resulting circuits can solve matrix equations in a single step in continuous time.

However, practical deployment remains challenging due to reliability and robustness issues, including device variability, resistance drift, limited endurance, and noise.\cite{8714954,9116244,9073995} These non-idealities can reduce computational accuracy and long-term system stability. Another major challenge is the reprogramming overhead required to represent different weights across neural network layers and models. Frequent weight updates lead to energy intensive write operations and increased latency, while repeated programming can accelerate device degradation because of the limited endurance of RRAM devices. In addition, accurately programming analog conductance states \cite{9474084,basisN} remains difficult, owing to nonlinear and asymmetric switching dynamics.
\begin{figure*}
    \centering
    \definecolor{solverWire}{HTML}{263238}
\definecolor{solverFeedback}{HTML}{2F6F78}
\definecolor{solverArray}{HTML}{ededed}
\definecolor{solverCell}{HTML}{5778a4}
\definecolor{solverCellEdge}{HTML}{5778a4}
\definecolor{solverMuted}{HTML}{5F6B70}

\tikzset{
  solver wire/.style={draw=solverWire, line width=0.85pt, line cap=round, line join=round},
  solver feedback/.style={draw=solverWire, line width=1.05pt, line cap=round, line join=round},
  solver junction/.style={circle, fill=solverWire, draw=none, inner sep=1.35pt},
  solver port/.style={circle, draw=solverWire, fill=white, line width=0.8pt, minimum size=4.2pt, inner sep=0pt},
  solver title/.style={font=\sffamily\small\bfseries, text=solverWire},
  solver sublabel/.style={font=\sffamily\scriptsize, text=solverMuted},
  pics/solver resistor/.style={code={
    \draw[solver wire, fill=white] (-0.34,-0.105) rectangle (0.34,0.105);
    \coordinate (-left)  at (-0.52,0);
    \coordinate (-right) at ( 0.52,0);
    \draw[solver wire] (-0.52,0)--(-0.34,0) (0.34,0)--(0.52,0);
  }},
  pics/solver opamp/.style={code={
    \draw[solver wire, fill=white] (-0.48,-0.43) -- (-0.48,0.43) -- (0.54,0) -- cycle;
    \coordinate (-inm) at (-0.48,0.18);
    \coordinate (-inp) at (-0.48,-0.18);
    \coordinate (-out) at (0.54,0);
    \node[font=\sffamily\scriptsize, text=solverWire] at (-0.31,0.18) {$-$};
    \node[font=\sffamily\scriptsize, text=solverWire] at (-0.31,-0.18) {$+$};
  }}
}

\begin{tikzpicture}[x=1cm,y=1cm]
  \def\arrayL{0.00}
  \def\arrayR{4.20}
  \def\arrayT{1.78}
  \def\arrayB{-1.78}
  \def\opX{5.72}
  \def\outX{8.48}

  \coordinate (c1) at (0.72,0);
  \coordinate (c2) at (2.10,0);
  \coordinate (c3) at (3.48,0);
  \coordinate (r1) at (0,1.18);
  \coordinate (r2) at (0,0);
  \coordinate (r3) at (0,-1.18);

  \node[solver title] at (2.10,2.4) {cross-point array $ A\in\mathbb{R}^{n\times n}$};
  \node[solver sublabel] at (2.10,2.05) {programmable conductances};
  \node[solver title] at (6.05,2.4) {operational amplifiers};

  \path[fill=solverArray, rounded corners=3pt]
       (\arrayL,\arrayB) rectangle (\arrayR,\arrayT);

  \foreach \y in {1.18,0,-1.18}{
    \draw[solver wire] (\arrayL-0.15,\y) -- (4.92,\y);
  }
  \foreach \x in {0.72,2.10,3.48}{
    \draw[solver wire] (\x,\arrayT+0.10) -- (\x,\arrayB);
  }

  \foreach \x in {0.72,2.10,3.48}{
    \foreach \y in {1.18,0,-1.18}{
      \draw[fill=solverCell, draw=solverCellEdge, line width=0.75pt]
        (\x,\y) circle[radius=0.205];
    }
  }

  \foreach \name/\y in {1/1.18,2/0,3/-1.18}{
    \node[solver port] (in\name) at (-2.05,\y) {};
    \pic (rin\name) at (-1.24,\y) {solver resistor};
    \draw[solver wire] (in\name) -- (rin\name-left);
    \draw[solver wire] (rin\name-right) -- (\arrayL-0.15,\y);
  }
  %\node[solver sublabel] at (-1.24,1.46) {$G_0$};

  \draw[decorate, decoration={brace,amplitude=4.5pt,mirror}, draw=solverMuted, line width=0.7pt]
    (-2.32,1.33) -- (-2.32,-1.33)
    node[midway,left=8pt] {$- b$};

  \foreach \name/\y in {1/1.18,2/0,3/-1.18}{
    \pic (oa\name) at (\opX,\y) {solver opamp};
    \draw[solver wire] (4.92,\y) -- (5.08,\y) -- (5.08,\y+0.18) -- (oa\name-inm);
    \draw[solver wire] (oa\name-inp) -- ++(-0.14,0) -- ++(0,-0.24);
    \draw[solver wire] ($(oa\name-inp)+(-0.25,-0.24)$) -- ($(oa\name-inp)+(-0.03,-0.24)$);
    \draw[solver wire] ($(oa\name-inp)+(-0.21,-0.30)$) -- ($(oa\name-inp)+(-0.07,-0.30)$);
    \draw[solver wire] ($(oa\name-inp)+(-0.17,-0.36)$) -- ($(oa\name-inp)+(-0.11,-0.36)$);
  }

  \draw[solver feedback]
    (7.94,1.18) -- (7.94,-3.12) -- (0.72,-3.12) -- (0.72,\arrayB);
  \draw[solver feedback]
    (7.66,0) -- (7.66,-2.74) -- (2.10,-2.74) -- (2.10,\arrayB);
  \draw[solver feedback]
    (7.38,-1.18) -- (7.38,-2.36) -- (3.48,-2.36) -- (3.48,\arrayB);

  \foreach \name/\y/\tapx in {1/1.18/7.94,2/0/7.66,3/-1.18/7.38}{
    \draw[solver wire] (oa\name-out) -- (\outX,\y);
    \node[solver junction] at (\tapx,\y) {};
    \node[solver port] (out\name) at (\outX,\y) {};
  }

  \draw[decorate, decoration={brace,amplitude=4.5pt}, draw=solverMuted, line width=0.7pt]
    (8.75,1.33) -- (8.75,-1.33)
    node[midway,right=8pt] {$x$};

  \node[solver sublabel, text=solverWire, fill=white, inner xsep=3pt]
    at (4.45,-3.12) {global negative feedback};

\end{tikzpicture}
    \caption{Analog crossbar circuit for solving the linear system $Ax = b$ using operational-amplifier feedback. The matrix $A$ is encoded in the cross-point conductances. The input vector $-b$ is applied as voltages across input resistors of unit conductance (used for matrix scaling) at the row terminals. The row lines of the crossbar are connected to the inverting inputs of operational amplifiers, which are held at virtual ground at steady state by negative feedback. The amplifier outputs are fed back to the column lines of the crossbar, and the steady-state amplifier output voltages represent the solution vector $x$. At steady state, Kirchhoff’s current law at each row summing node requires the net current from the crossbar and the input branch to be zero (\ie, $Ax-b=0$), yielding $x=A^{-1}b$ for invertible $A$.
}
    \label{fig:linear_system}
\end{figure*}
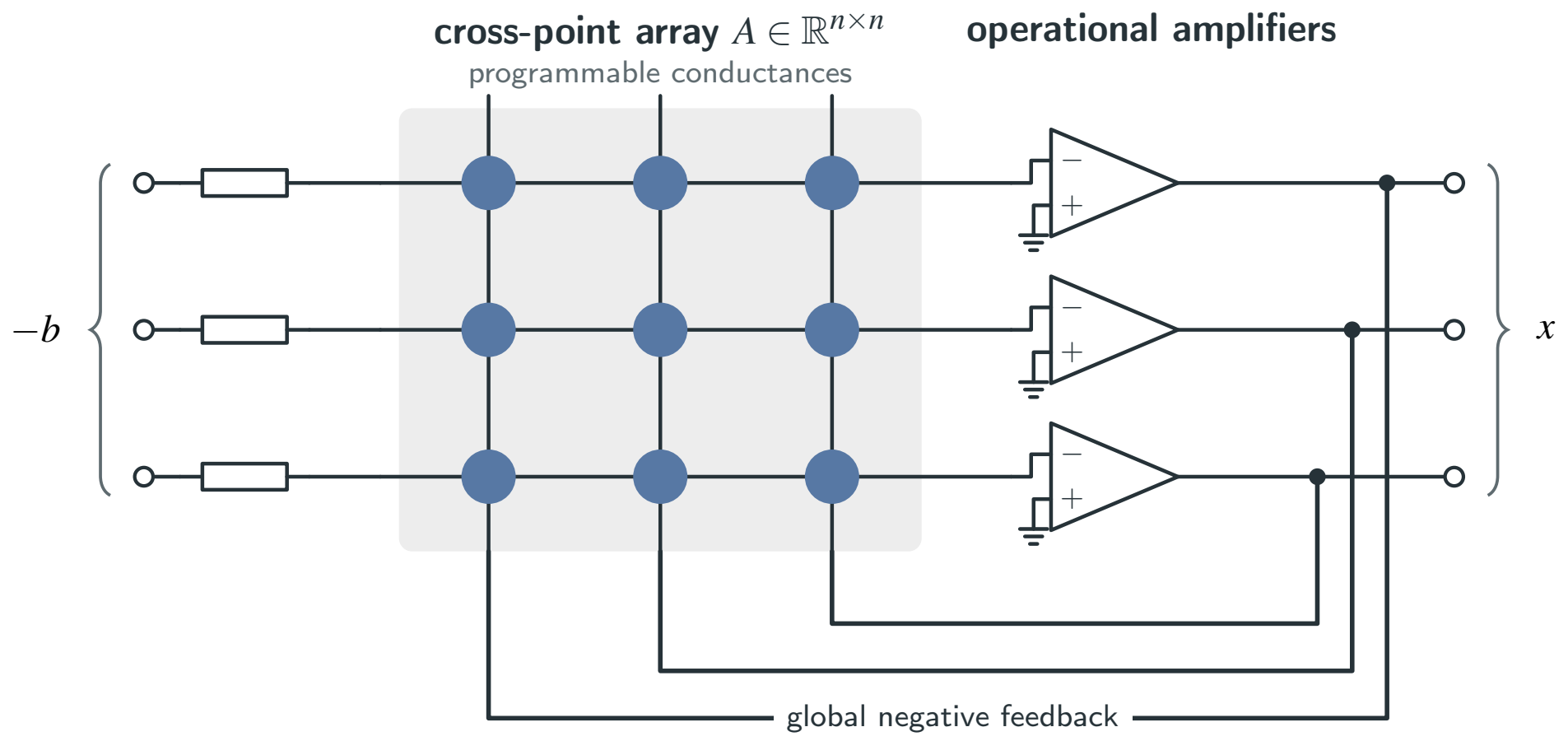

Beyond MVMs, in-memory computing can directly solve matrix equations such as linear systems, matrix inversion, generalized inverses, and matrix factorization. Two recent review papers\cite{sun2022invited,miao2026exploiting} have comprehensively surveyed in-memory analog circuits for these primitives, showing that closed-loop feedback configurations of memristive crossbars and operational amplifiers can compute solutions in constant time ($\mathcal{O}(1)$ complexity) by physically implementing iterative algorithms as continuous-time dynamics (see Fig.~\ref{fig:linear_system}). More recently, this paradigm has been extended to in-memory non-negative matrix factorization.\cite{wang2026memory} By mapping alternating non-negative least squares subproblems onto memristive crossbar arrays, each subproblem is solved in one step, which is useful in decomposition tasks arising in image compression and recommender systems. Analog matrix inversion circuits have also been employed to accelerate second-order training of neural networks,\cite{zhang2025first} where the crossbar-based inverter computes a preconditioner or inverse-Hessian-vector product, significantly reducing training time compared to pure digital methods.

In-memory matrix equation solvers share the reliability challenges of MVM accelerators, including device variability, programming inaccuracy, and analog-to-digital converter (ADC) / digital-to-analog converter (DAC) overhead, but they also face additional challenges regarding precision and scalability. For ill-conditioned problems (\ie, problems with a high condition number), even small conductance errors are amplified. To enhance precision, bit slicing and iterative refinement, such as Richardson or Newton iterations, can be adopted.\cite{zuo2025precise} In these approaches, memory devices are programmed to a small number of distinguishable conductance states with sufficient readout margins (\eg, 4-bit precision).\cite{li2024maximum} This combination reduces both variability-induced error and condition-number sensitivity, although ill-conditioned matrices may still require preconditioning techniques to improve the effective condition number.\cite{feinberg2021analog} Scalability is another issue because practical crossbar dimensions are limited, and handling large matrices necessitates partitioning. One potential partitioning approach is the BlockAMC scheme,\cite{pan2024blockamc} which decomposes a high-dimensional matrix into smaller blocks that sequentially undergo MVM and matrix equation solving to recover the full solution vector. Combining high-precision bit-slicing with BlockAMC enables large-scale, high-accuracy matrix equation solving for real-world applications.

ADC precision is a further design concern for matrix equation solvers. A detailed analysis of MVM precision\cite{li2025constraints} shows that the output effective number of bits is limited by device conductance precision, input DAC resolution, and the number of simultaneously activated rows. Consequently, the required ADC resolution is modest. For matrix equation solving, high precision is achieved through bit slicing and iterative refinement.\cite{zuo2025precise} Decomposing the matrix into low-precision slices restricts the MVM output dynamic range, and the closed-loop analog solver inherently produces an approximate result, further relaxing ADC precision demands. Therefore, a modest ADC resolution (\eg, 4--8 bits) suffices, allowing designers to prioritize area, speed, and power efficiency. With appropriate cross-layer strategies, ADC overhead can thus be effectively mitigated in high-precision in-memory matrix equation solvers.
\subsubsection{Machine learning}
In machine learning applications, training is usually performed infrequently, but inference is repeated for every query and, at scale, can account for more than 90\% of the total energy and cost.\cite{desislavov_trends_2023} This makes inference the primary target for in-memory computing. Mapping a software-trained network onto analog hardware, however, is not a lossless operation. The conductances that store the weights are subject to device-to-device variability as well as nonlinear and imperfect programming, so the weights actually realized in the crossbars deviate from their intended values and the accuracy of the network degrades.\cite{aguirre_hardware_2024} Thus, a network cannot simply be trained once and have its weights copied from one hardware instance to another.\cite{hinton_forward-forward_2022} This regime has been termed mortal computation, as the learned solution is bound to the properties of the physical substrate and does not outlive it. The model must therefore be trained, or at least fine-tuned, on the specific system that will run it. Even then, the weights may drift over time during operation, in which case periodic retraining or recalibration becomes necessary.\cite{chang_mitigating_2021}

At the system level, multilayer neural-network accelerators face another challenge, namely the ``ADC wall'' in which the energy and area required for analog-to-digital conversion can exceed the cost of the analog computation itself.\cite{Shafiee2016ISAAC,Ankit2019PUMA,yao2020fully,RACA2024} An ideal analog accelerator would therefore keep intermediate activations entirely in the analog domain and perform digitization only at the final network output. However, realizing such architectures requires analog representations to remain robust against device variations and noise across multiple layers. This motivates the development of analog-domain error mitigation and self-calibration techniques, including on-device continual learning and adaptive calibration mechanisms that compensate hardware non-idealities without frequent conversion to the digital domain.

During training, backpropagation, the standard method for computing gradients, is poorly suited for analog hardware.\cite{kendall_training_2020} In the analog domain, backpropagation usually requires dedicated circuitry, and the device mismatches and nonidealities present in such hardware cause gradient errors. This can be mitigated by computing the activations and their derivatives in the digital domain, but doing so requires converters between every pair of layers and thus reintroduces the ADC bottleneck discussed above. Equilibrium propagation\cite{scellierEquilibriumPropagationBridging2017} and coupled learning\cite{stern_supervised_2021}, two closely related approaches for energy-based models, avoid these problems by treating the network as a dynamical system that settles into an equilibrium state of minimum energy for a given input. Training runs the system twice, both times using the same circuit and the same dynamics. In the first phase, it relaxes with its outputs left free. In the second phase, the outputs are driven toward the desired targets by a weak nudging force in equilibrium propagation and by clamping them to slightly shifted values in coupled learning. Each weight is then updated from the difference of the squared voltage drops measured across its device in the two equilibrium states.\cite{kendall_training_2020} Since these updates are purely local, no error signals have to be transported through the network, and device nonidealities do not prevent reliable gradient computation.
\subsection{\label{sec:bio_computing}Biological computing---By Frank Delvigne, Alex Fedorec, and Daniel Garza}
Computation is a defining feature of living systems, which are often characterized by their complexity,\cite{doi:10.1073/pnas.1807890115,mitchell2009complexity} far-from-equilibrium dynamics,\cite{cao2025stochastic} use of free energy,\cite{schrodinger1944what} and capacity for Darwinian evolution.\cite{benner2010defining} Although computation is not unique to life, it may be one of life's most general properties. If multi-scale learning and information processing are unifying features of life,\cite{hogeweg2011roots,vanchurin2022toward} then computation (\ie, the transformation of inputs and internal states into functional outputs) is a basic concept they share. 

The properties of analog computers discussed in Sec.~\ref{sec:intro} have close biological counterparts. For example, living systems use special-purpose components that evolved to solve specific tasks. These components can further evolve, replicate, and form larger computing structures. We therefore argue that biological systems form a distinct class of analog computers.
\subsubsection{A minimal model for biological computing}
To develop a unified theory of biological computation, one possibility is linking replication, transcription, and translation with metabolism,\cite{costa2021remodel} evolution, and ecology.\cite{takeuchi2025generalizing} For analog computing, however, a minimal working definition is already useful. At a minimum, biological computation occurs when a living system uses inputs such as temperature, chemical gradients, pressure, or combinations of these to produce outputs such as growth, movement, secretion, differentiation, death, or reproduction.

The corresponding minimal model assumes a functional relation in which input and state variables constrain the outputs of the considered biological system. Let $X(t)$ denote the input variables that act on the system at time $t$. The quantity $S(t)$ contains the state variables that condition the input-output relationship, while $Y(t)$ contains output variables. Because biological systems may be stochastic, this relation can then be written as $Y(t)\sim P_{\theta(t)}(X(t),S(t))$, where $P_{\theta(t)}$ is the probability distribution over possible outputs for the current inputs and states, and $\theta(t)$ contains the parameters and structural constraints that determine how inputs and states are mapped to outputs. Mutual information provides one measure of the statistical dependence between inputs and outputs,\cite{tkavcik2016information} although functional coupling also depends on mechanism, state, and context. 

In biological computation, we distinguish between weak computation and strong computation. In weak computation, the model parameters $\theta(t)$ are kept constant over the timescale of the computation. In strong computation, the model parameters evolve according to some function $G(\cdot)$. That is, $\theta(t+\Delta t)=G(\theta(t),X(t),S(t),Y(t))$. The underlying computational unit is not centralized but a network of coupled chemical, physical, spatial, ecological, or biological variables whose interactions determine both the computation and adaptation.

At evolutionary timescales, the update function $G$ captures changes in $\theta(t)$ across generations. When the biological organization that implements the input--state--output mapping is reproduced with variation, and its variants differ in persistence or reproductive success, that organization also constitutes a unit of selection.\cite{lewontin1970}  Under these conditions, the computational and evolutionary units coincide. This connection links biological computation directly to evolutionary dynamics, providing a framework that can accommodate multi-level selection,\cite{10.1111/jeb.12566} complexification through time,\cite{aguerayarcas2025what} and major evolutionary transitions\cite{szathmary1995major,vanchurin2022toward} while remaining compatible with neutral and nearly neutral variation.\cite{kimura1968evolutionary}

The described framework of biological computation applies to biological systems across scales, including molecules,\cite{bray1995protein} cells,\cite{perkins2009strategies} populations,\cite{cornforth2014combinatorial,scialdone2013arabidopsis} communities,\cite{Ostovar_2020,karkaria2020microbial} ecosystems,\cite{ushio2023computational} and the biosphere.\cite{falkowski2008microbial,sole2015bioengineering}
\subsubsection{Measuring biological input--output relationships}
\begin{figure}
    \centering
    \includegraphics[width=\linewidth]{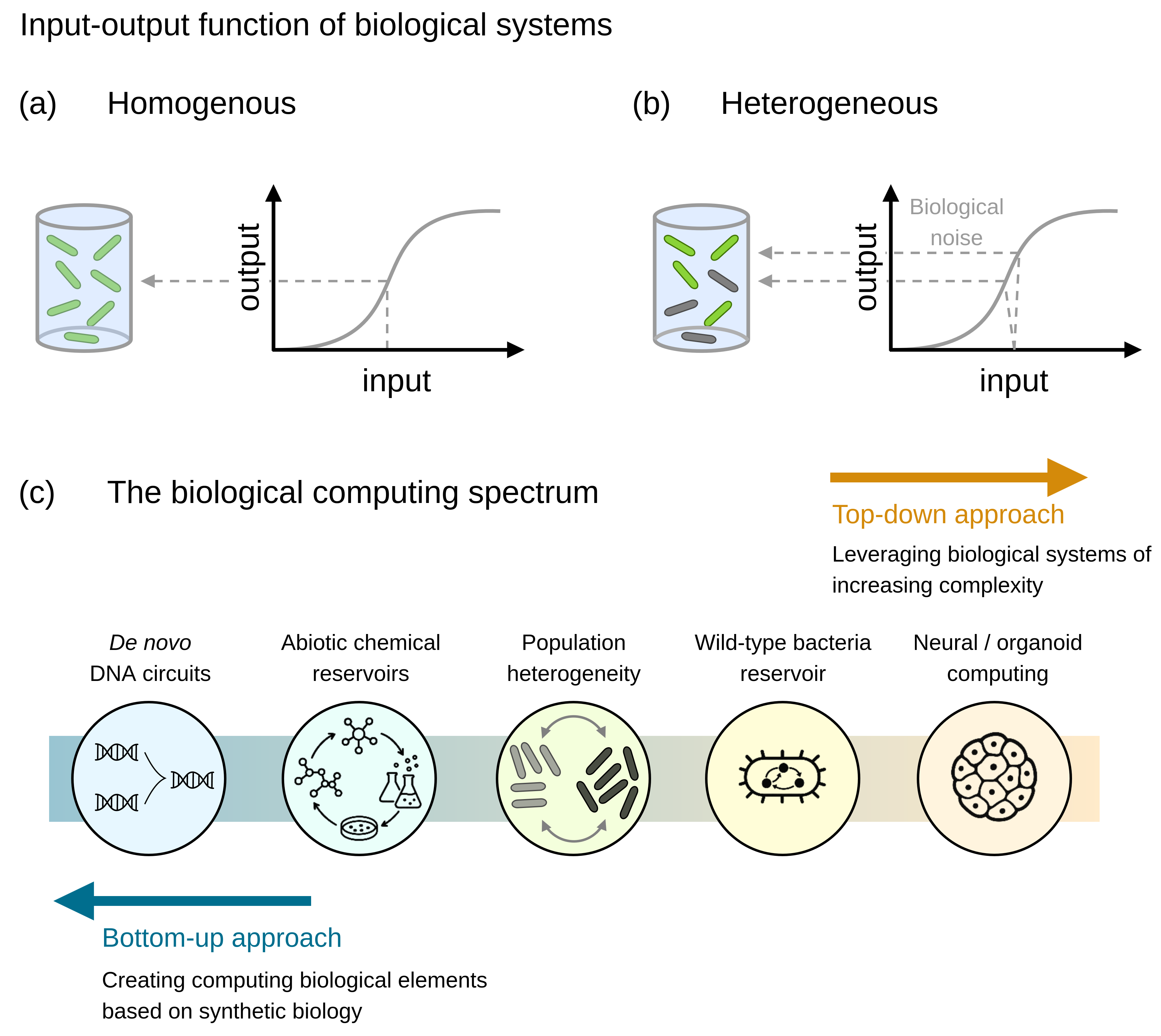}
    \caption{Schematic of biological computing approaches. (a) A single input signal applied to a homogeneous microbial population produces a deterministic output response. (b) Biological noise in heterogeneous populations generates a distribution of outputs across subpopulations. (c) Biological computing approaches range from fully engineered de novo DNA circuits (bottom-up, synthetic biology) to biological substrates of increasing complexity (top-down), spanning abiotic chemical reservoirs, population heterogeneity-based reservoirs, wild-type bacterial reservoirs, and neural/organoid computing systems.}
    \label{fig:bio_computing}
\end{figure}

To empirically examine input-output relationships of specific biological functions, one has to carry out measurements that usually rely on natural or synthetic gene circuits. In the synthetic context, fluorescent protein expression is the most widely used readout, with common inputs including nutrient changes, light stimulation via optogenetic systems, or the addition of inducers. There are an abundance of other outputs that one can track, including in natural biological systems, such as metabolic responses, proliferation, motility, aggregation, or ecological interactions. In homogeneous populations, a given input maps to a single deterministic output [see Fig.~\ref{fig:bio_computing}(a)]. However, biological noise, arising from stochastic gene expression and phenotypic heterogeneity, causes isogenic populations to produce a range of outputs in response to the same input [see Fig.~\ref{fig:bio_computing}(b)]. This property of biological systems is both a challenge for classical measurement frameworks and, as argued here, a latent computational resource.
\subsubsection{Biological computation as a source of analog computing architectures}
Biological systems perform analog computation based on information encoded in continuous variables such as molecular concentration,\cite{su2023bacterial} signal duration,\cite{dessaud2007interpretation} oscillation frequency,\cite{zhang2016processing} and spatial patterns.\cite{dubuis2013positional} Information processing occurs via chemical reactions, diffusion, spatial organization, and stochastic fluctuations. Features often viewed as constraints in engineered systems, including noise, transport limitations, and local interactions, can be computational resources in biological systems.\cite{eldar2010functional,kondo2010reaction}

In tissues and biofilms, cell positions are constant or change slowly relative to their neighbors. Diffusing molecules form gradients and boundaries that shape which cells interact and how strongly.\cite{backer2026spatial} Because each cell responds from its own position and history, spatial organization can encode memory in position-dependent states and enable distributed processing through local interactions. For instance, in \textit{Bacillus subtilis} biofilms, localized optical stimuli trigger potassium-mediated changes in membrane potential. These changes affect how the stimulated cells synchronize with global ion oscillations. The stimulated region remains distinguishable for more than three hours, which shows that a transient spatial input can be converted into persistent spatially-encoded memory.\cite{yang2020encoding} 

In biological systems, noise can contribute directly to information processing.\cite{malaguti2021theory} For example, during bacterial chemotaxis (\ie, movement in response to chemical signals), a cell explores its environment through periods of unidirectional motion and random changes in direction. The signaling pathway uses this randomness by lengthening runs when conditions improve and shortening them when conditions worsen.\cite{mattingly2021escherichia} Random reorientation contributes to exploration, while signal-dependent changes in run length contribute to exploitation. Another example illustrating the role of noise in biology is regulatory switching, where noise can move a system between alternative phenotypes.\cite{ham2025mapping} In microbial populations, stochastic switching produces heterogeneity that supports bet-hedging, creating subsets of cells that may already occupy a phenotype suited to a future environment.\cite{kussell2005phenotypic}

One can also implement elementary mathematical operations using biochemical processes.\cite{buisman2009computing,sarpeshkar2014analog} For example, a Michaelis--Menten process produces a saturating input-output response. If two genes independently produce the same output protein, their production rates add. Under first-order degradation, the steady-state protein concentration is the sum of the concentrations each gene would produce alone, giving a molecular form of addition. Under mass-action kinetics, the rate of molecular association is proportional to the product of the reactant concentrations, giving a molecular form of multiplication. More generally, nonlinear reaction kinetics can generate higher-order operations, including square-like and square-root-like relationships under appropriate conditions.\cite{anderson2025chemical}
\subsubsection{Biological systems as a substrate for analog computation}
Biological systems, including biomolecules and individual cells, multicellular populations, and tissues, provide a substrate for implementing analog computation. In contrast to conventional electronic hardware, biological computational substrates are active, adaptive, and self-organizing. 

At the molecular level, information is represented by continuous variables such as molecular concentrations, binding affinities, and reaction rates. Corresponding analog computations can be implemented using engineered chemical reaction networks,\cite{araujo2023universal} molecular circuits,\cite{daniel2013synthetic,teo2015synthetic} and DNA strand-displacement.\cite{song2016analog}  For example, DNA strand-displacement circuits can represent numerical values as strand concentrations, and elementary reaction modules for addition, subtraction, and multiplication can be combined to perform more complex computations.\cite{song2016analog} Such systems represent a bottom-up approach to biological computing in which computational functions are specified through the design of individual molecular interactions.\cite{cherry2025supervised}

At the cellular level, one can use the natural properties of cells to employ them as sensors, processors, and controllers. For example, synthetic gene circuits have been engineered to perform classification, store memory, count environmental events, and produce other specific outputs.\cite{xie2018designing} Most effort in the engineering of cells to compute has revolved around producing digital logic,\cite{nielsen2016genetic} though there have been some notable attempts to produce to produce analog biocomputation.\cite{daniel2013synthetic,rubens2016synthetic}

At the multicellular level, engineered microbial systems can perform and distribute computational tasks across populations, which can sense inputs, transform chemical states, and communicate through metabolism, secretion, cross-feeding, and quorum sensing.\cite{kylilis2018tools} While a well-mixed community can perform distributed computation through a globally accessible environmental biochemical ``memory'', the addition of spatial organization can determine which interactions occur and therefore becomes part of the computational architecture. In such systems, computation emerges from interactions among cells. Thus, the population as a whole becomes the computing unit.\cite{regot2011distributed} In one design, receiver colonies with simple activation functions, placed at precise distances from diffusible morphogen sources, encode distinct logic gates, allowing the same molecular signal to compute different functions depending on the geometry. Simply by changing the spatial pattern of the inputs and outputs, one changes the logic function encoded by the system.\cite{fedorec2024emergent} Computation therefore emerges not only from molecular identity but also from spatial arrangement. Other examples of multicellular computing substrates include organoids, engineered tissues, and neural cultures.\cite{robbins2026goal} 

One distinctive feature of biological computing is that the substrate is itself alive. However, the adaptive properties that make biological systems an attractive computational substrate can complicate their engineering. For this reason, biological computing should be considered a specialized computational paradigm for applications in living systems. Near-term applications include biosensing, bioproduction, bioremediation, therapeutic control, and living materials, where cells operate in the same environment that defines the computational problem.\cite{grozinger2019pathways}

In Fig.~\ref{fig:bio_computing}(c), we classify biological computing approaches according to their level of biological organization. At one extreme, \textit{de novo} DNA circuits implement computation through molecular interactions in cell-free systems.\cite{cherry2025supervised} Intermediate approaches include abiotic chemical reservoirs such as the formose reaction, which are based on self-organizing dynamics without circuit-level design.\cite{baltussen2024chemical} Another example involving an intermediate degree of biological organization is wild-type bacterial reservoir computers, which repurpose unmodified organisms whose intrinsic metabolic dynamics provide a reservoir.\cite{AHAVI2026101654} At the other extreme are wetware and organoid intelligence. These systems are based on intact biological neural tissue, where computation emerges from electrophysiological activity and adaptive network reconfiguration.\cite{LiBrainOnChip25} 

This range of biological computing systems reflects two complementary design principles. Bottom-up approaches build computing elements from synthetic parts. Such approaches offer precision and mechanistic understanding at the cost of design complexity and limited scalability. Top-down approaches leverage the natural dynamics of increasing biological complexity, trading molecular control for emergent computational richness.
\subsection{\label{sec:general_purpose}General purpose analog computing---By Lucas Wetzel, Sara Obergassel, and Pavol Pitonak}
Traditionally, there were two types of analog computers: (i) special-purpose machines hardwired for one task, which made up the bulk of units ever built; and (ii) general-purpose machines reprogrammed by manually connecting computing elements by means of patch cables for each new problem. Either way, configuration was physical. No software layer separated the problem from the hardware.\cite{ulmann2026analog} The vision behind a general-purpose analog computer (GPAC) is different: a single reconfigurable substrate that can be programmed to represent a wide class of dynamical systems, much like a digital processor is programmed to run a wide class of algorithms. The relevant elementary unit is not a fixed circuit, but a set of analog building blocks together with a configuration that connects them.

Today, the vision of developing a GPAC is pursued through, \eg, hybrid analog-digital architectures.\cite{koppel2021} The analog part integrates the continuous dynamics and carries out the computation. The digital part configures it, sets coefficients, supervises timing, reads out results, and closes control loops that the analog part cannot close on its own. Neither part is sufficient alone, and the hybrid, not the analog core in isolation, is the object of interest.
\subsubsection{Building blocks of a hardware GPAC}
\label{sec:gpac_building_blocks}
The theoretical model of a GPAC is described in Sec.~\ref{sec:new_quest}. It consists of constant units, summers, multipliers, and integrators. The corresponding wiring fixes a polynomial vector field. Additionally, comparators and arbitrary function generators are provided. Delay lines, convolution units, and nonlinear activation blocks, such as the rectified linear unit and a hard hyperbolic tangent, are currently under test and are expected to become available on the planned GPAC platform.

Comparators and activation blocks bring thresholding and switching into the analog domain. Delay lines and convolution units are useful to account for signal history and memory effects. In the GPAC design we pursue, building blocks are organized into interchangeable units connected in a reconfigurable manner, so a larger problem is handled by adding units rather than by redesigning the machine. This modular construction is what makes the same hardware general-purpose: the configuration changes, the hardware does not.
\subsubsection{The REDAC platform}
\label{sec:redac_platform}
\begin{figure}
    \centering
    \includegraphics[width=\linewidth]{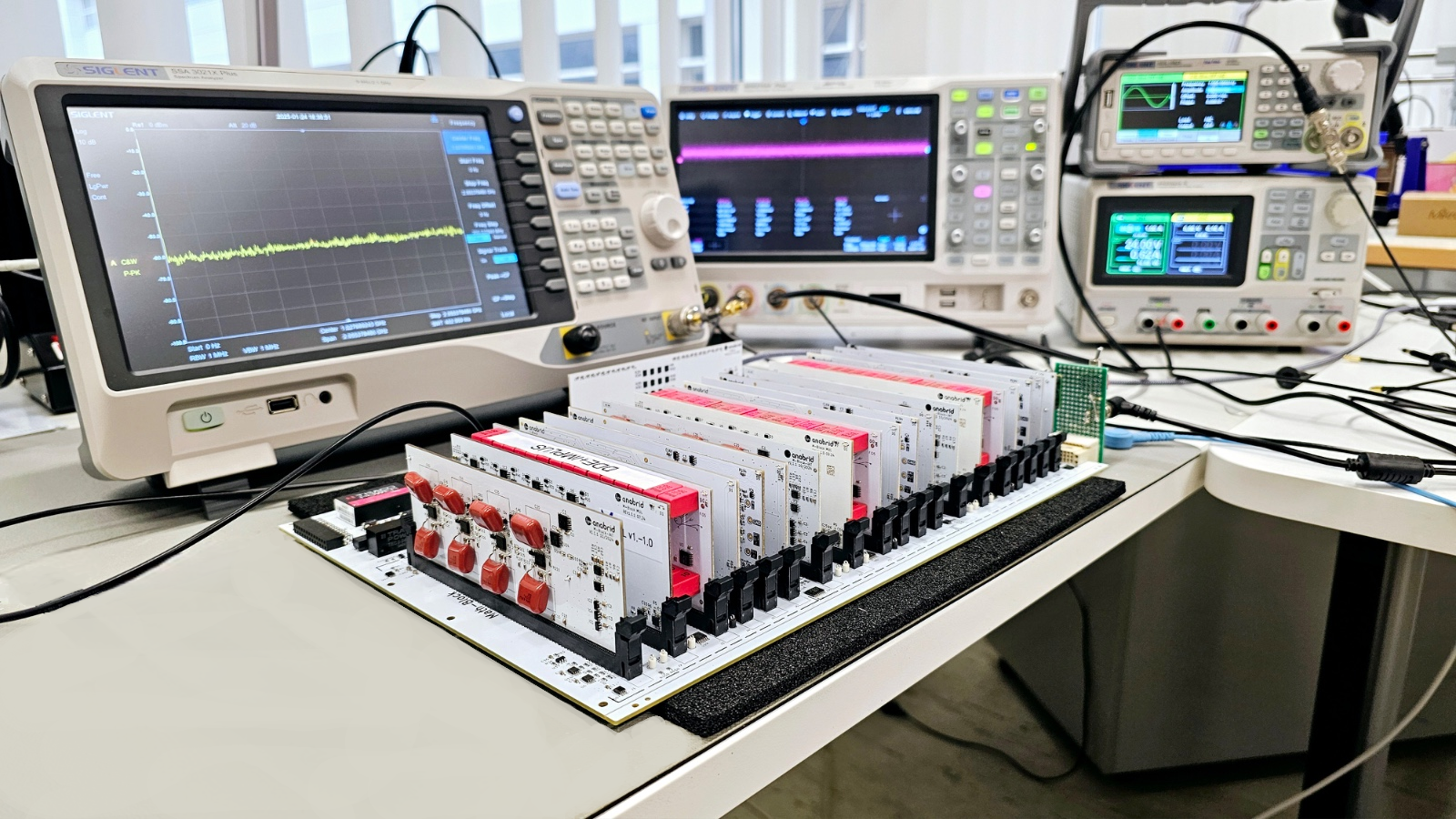}
    \caption{One mREDAC module of the reconfigurable discrete analog computer (REDAC). The individual groups of computing and routing elements are implemented using dual inline memory module (DIMM) form factor elements. A module like this typically contains 24 integrators, 12 multipliers, 96 digitally controlled coefficents, and the necessary switch matrices for automatic reconfiguration under control of a digital computer. Figure courtesy of anabrid GmbH.}
    \label{fig:redac}
\end{figure}
REDAC is a reconfigurable discrete analog computer (see Fig.~\ref{fig:redac}). It operates in continuous time and has continuous variables. REDAC is a development platform and built to study hybrid systems in a way that integrated hardware does not yet offer. Because every building block is accessible and measurable, REDAC enables developers to examine accuracy, bandwidth, and energy use under controlled conditions. It is also a useful platform for testing programming models and compiler protocols before they are committed to silicon.

\begin{figure*}
    \includegraphics[width=0.9\textwidth]{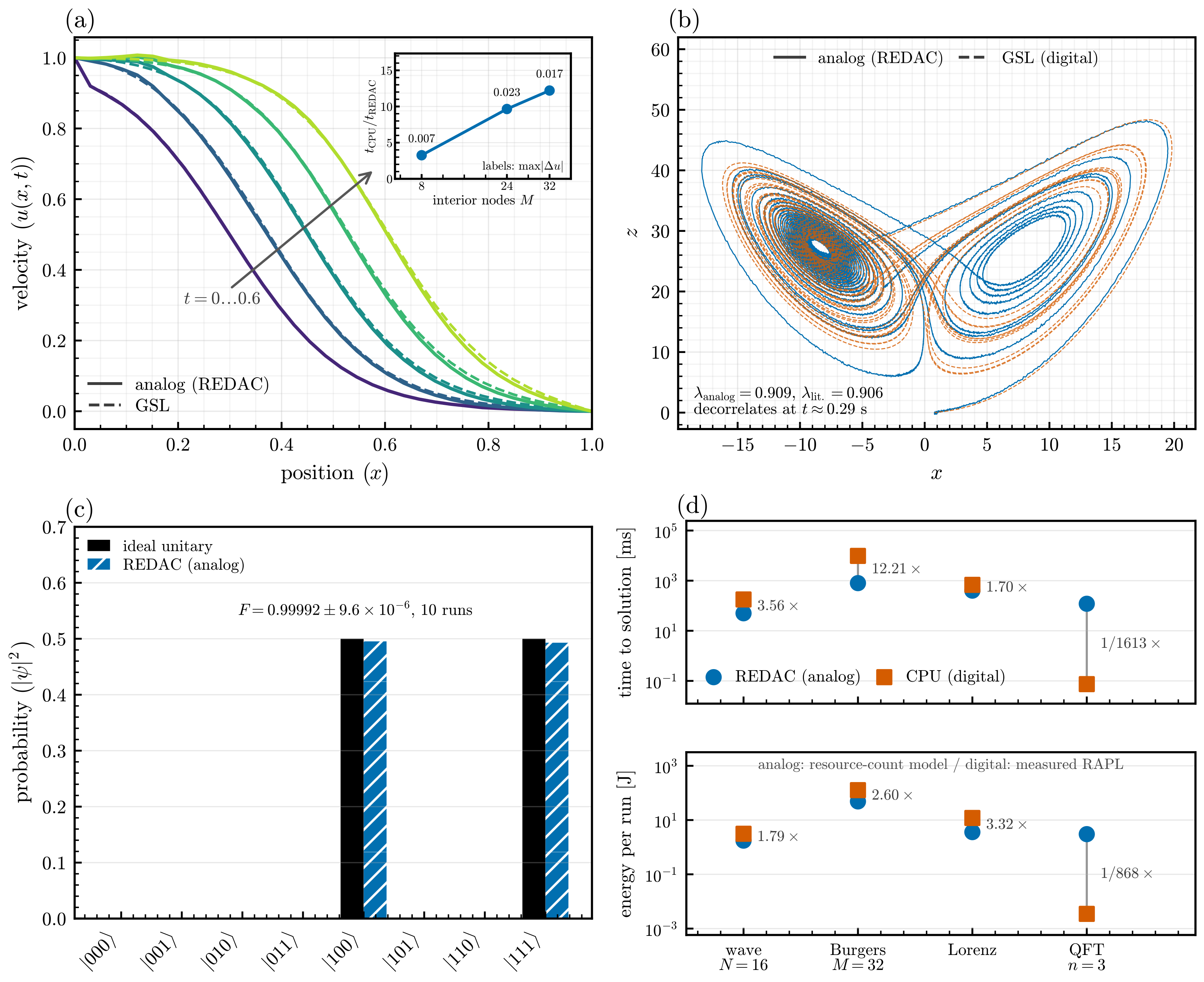}
    \caption{Representative computations on the REDAC platform across application classes, each compared against a digital reference solution. All analog data were measured on the physical device.
    (a) Viscous Burgers equation with a traveling shock. The analog implementation uses $M=32$ interior nodes across four carriers and one hardware integrator per node. We show analog (solid) and the GSL reference (dashed) velocity profiles at five times. The maximum deviation over the full window is $3.8\,\%$ of the unit velocity jump. The inset sweeps the mesh over $M=8$, $24$ and $32$. The speed advantage over the processor increases with $M$, from a factor $3.3$ to a factor $10.3$.
    (b) Lorenz attractor ($\sigma=10$, $\rho=28$, $\beta=8/3$) in the $x$--$z$ plane. The Lyapunov exponent fitted to the analog run is $0.909$. The corresponding literature value is $0.906$. Analog and digital trajectories decorrelate at $t\approx0.29$~s, as expected for a chaotic system.
    (c) Gate-based quantum Fourier transform on three qubits, applied to one full period of a sinusoid across the eight basis states. We show the output-state probabilities against the exact unitary, with mean fidelity $F=0.99992\pm9.6\times10^{-6}$ over ten realizations. The transform of that input is a line spectrum of weight $1/2$ at $\ket{100}$ and $\ket{111}$, and the six states that should be empty carry a total leakage below $10^{-4}$.
    (d) Time to solution and energy per run for all four cases. The annotation gives the digital-to-analog ratio. Analog time to solution is the chip op-time and analog energy is a resource-count model, whereas the digital figures are measured wall time and measured RAPL energy. The two energy figures are therefore not like-for-like. The Burgers and quantum Fourier transform op-times were increased to achieve higher accuracy, using coefficient pre-division and a slower integration rate. Therefore, the energy of those two cases is set by the accuracy target rather than by a floor of the hardware. The digital part of the transform is the exact unitary cascade, not an integrated solver, and costs $3.5$~mJ; the same circuit costs $11$~mJ on the \texttt{qiskit} state-vector simulator. The three-qubit transform sits far
    below the problem size at which the analog route pays off, and is shown to mark that boundary rather than to claim an analog advantage. Figure courtesy of anabrid GmbH.}
    \label{fig:redac_results}
\end{figure*}
Figure~\ref{fig:redac_results} shows four types of computations carried out on a REDAC, each measured against a digital reference. The reference for the three differential-equation cases is the GNU Scientific Library (GSL), integrated with a fixed fourth-order Runge--Kutta step set equal to the analog sampling interval, so that the two solutions are compared at the same temporal resolution. It was run on a 13th-generation Intel Core i7-13700H processor. The employed analog circuits are small. The Burgers case occupies a single carrier using $8$ integrators and $14$ multipliers, while the Lorenz case occupies a single carrier using only $3$ integrators and $2$ multipliers. The viscous Burgers shock is reproduced to within $0.7\,\%$ of the velocity jump over the full capture window. For the vibrating string, the deviation stays within $6\,\%$ of the amplitude over the first period. It grows over the twenty periods captured, as the analog and digital solutions accumulate a phase difference. For the Lorenz system, the Lyapunov exponent fitted to the analog run is $0.909$. This value is close to the value of $0.906$ reported in the literature. The two trajectories decorrelate at $t\approx0.29$~s, as expected for a chaotic system. A gate-based quantum Fourier transform on three qubits reaches a mean fidelity of $0.99948\pm4.1\times10^{-5}$ over ten realizations. Time to solution is the op-time of the analog circuit. It lies between $50$ and $400$~ms for these cases, a factor $1.7$ to $3.6$ below the measured processor wall time.

For the three differential-equation cases, the analog energy per run is a factor $1.8$ to $3.8$ below the digital processor energy. For the analog implementations, each carrier draws $7.8$~W, each integrator $0.115$~W, and each multiplier $0.42$~W. We compute the total energy per run by integrating these values over the op-time. The energy consumption of the digital Intel processor is measured by reading out the running average power limit (RAPL) counters,\cite{david2010rapl,khan2018rapl} using the pyJoules library.\cite{pyjoules} The analog and digital energy values are therefore not directly comparable, and we do not report their ratio as a measured energy ratio.

Two of the op-times in Fig.~\ref{fig:redac_results}(d) were increased to achieve higher accuracy. The Burgers coefficients are pre-divided in software, so that every emitted coefficient stays inside the representable range. The quantum gates are evolved at the slower of the two integration rates the hardware offers, which reduces the relative quantization error of the coefficient path. Under a fixed power model the energy scales with the op-time. The energy of these two cases is therefore set by the accuracy target, and is not a floor for the hardware. This trade between accuracy and energy is itself a design rule that a development platform is meant to expose.

The three-qubit transform is the counter-example, and we show it for that reason. It costs more time and more energy on the analog platform than on the processor. The crossover for this problem class therefore lies above three qubits.

These are single instances at small problem sizes, and we do not extrapolate them to a scaling claim. Applications of analog and hybrid platforms to dynamical systems, PDEs, and optimization are discussed in Sec.~\ref{sec:applications}.
\subsubsection{Toward integrated analog computing chips}
\begin{figure}
    \centering
    \includegraphics[width=\linewidth]{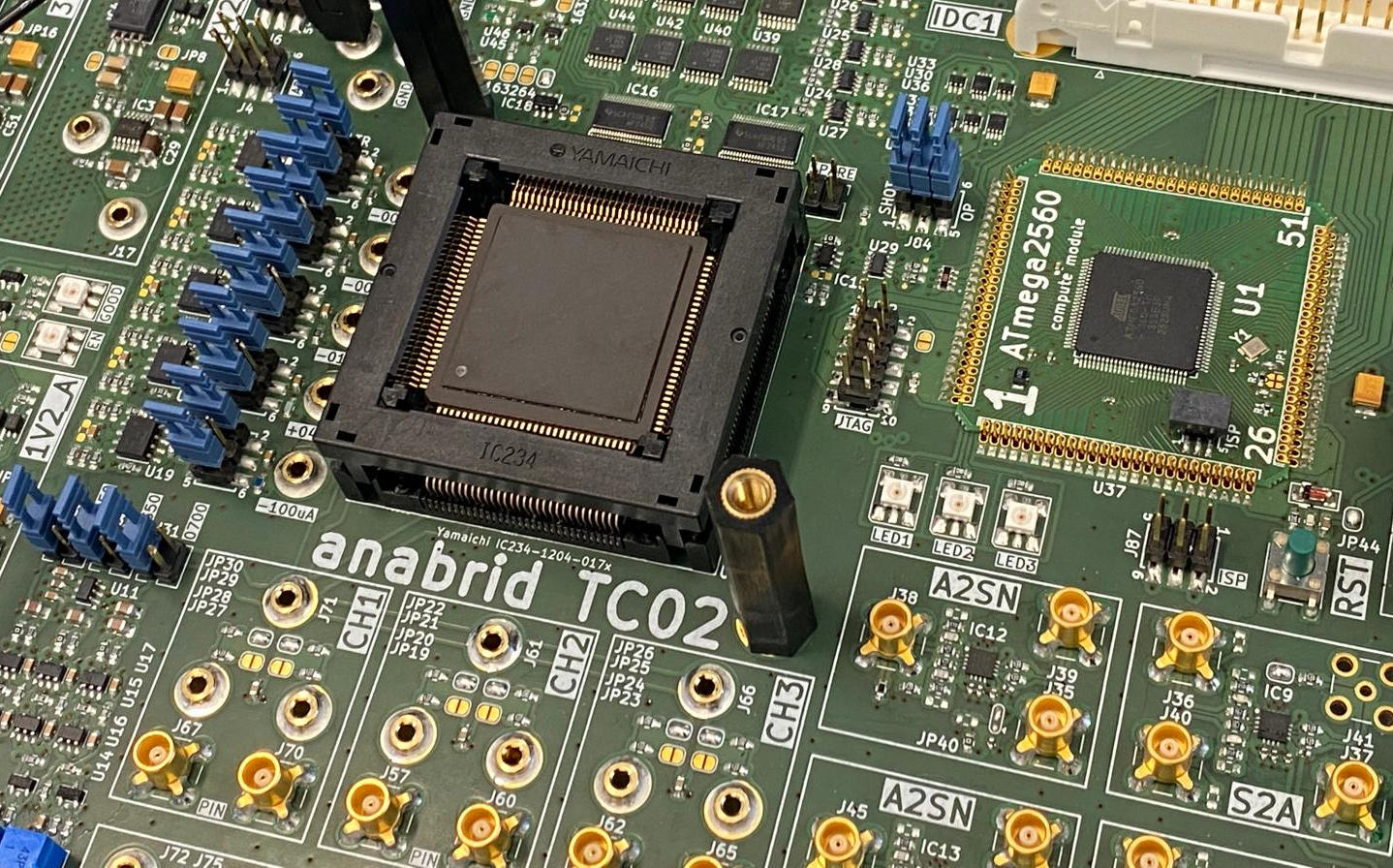}
    \caption{Prototype chip of a fully reconfigurable analog computer in 65 nm CMOS technology during characterization work in the lab. Figure courtesy of anabrid GmbH.}
    \label{fig:anabrid_chip}
\end{figure}

The current REDAC-type analog computing platforms are useful for benchmarking purposes. A natural next step is to develop compact and energy-efficient configurable analog computing systems on chips. However, rather than replicating large discrete analog computers on silicon, the objective is to develop a set of reusable analog and mixed-signal computational building blocks that can be configured to implement a wide variety of continuous-time dynamical systems. These building blocks include integrators, summing amplifiers, multipliers, comparators, programmable coefficient networks, function generators and configurable routing matrices. In conjunction with digital systems, the analog hardware can be configured by programming the interconnects, coefficients and operating parameters. In this way, one can realize different computational models on the same physical device without modifying the underlying circuitry.\cite{Killat2025,Co2006} Figure \ref{fig:anabrid_chip} shows the prototype of a fully reconfigurable analog computer on chip in 65 nm CMOS technology.

Implementation in complementary metal-oxide-semiconductor (CMOS) technology has two main advantages over discrete realizations. First, monolithic integration reduces parasitic capacitances and interconnect delays. These improvements enable lower power consumption, higher operating bandwidth, improved reliability, and greater computational density. Second, reusable analog building blocks are useful for systematically designing configurable analog and mixed-signal application-specific integrated circuits (ASICs).\cite{Gu2016}

However, integrated analog computing also faces significant design challenges. For example, device mismatch, process variations, thermal drift, noise, limited signal swing, and routing parasitics negatively impact accuracy and robustness. It is therefore important to complement integrated analog circuits with appropriate calibration and error-compensation methods to ensure reproducible operation and establish error bounds. These methods are especially relevant when scaling up the number of analog building blocks in integrated circuits.
\subsubsection{The wider landscape}
The REDAC platform is one of several current efforts to build reconfigurable or general-purpose physical computers. Another recent approach to developing programmable analog computers in standard CMOS is based on floating-gate cells.\cite{hasler2024} Additionally, coupled-oscillator circuits have been employed to implement Ising machines that solve combinatorial-optimization problems with reported speed and energy advantages over digital solvers.\cite{moy2022,lo2023} Physical reservoir computing is another related approach and uses a fixed nonlinear substrate, which can be electronic, photonic, or mechanical, with only a trained digital readout.\cite{tanaka2019} These systems differ in what they fix and what they expose to the programming layer. The distinguishing feature of the REDAC approach is the combination of a general extended-GPAC building-block set with a hybrid analog-digital architecture and a full software stack.

The software stack that turns REDAC-type hardware into a programmable instrument, from firmware and assembler through the compiler to cloud-based management, is described in Sec.~\ref{sec:software}.
\section{\label{sec:software}Software and systems development---By Christopher Jung and Daniel Thuerck}
\begin{figure*}
    \includegraphics[width=0.8\textwidth]{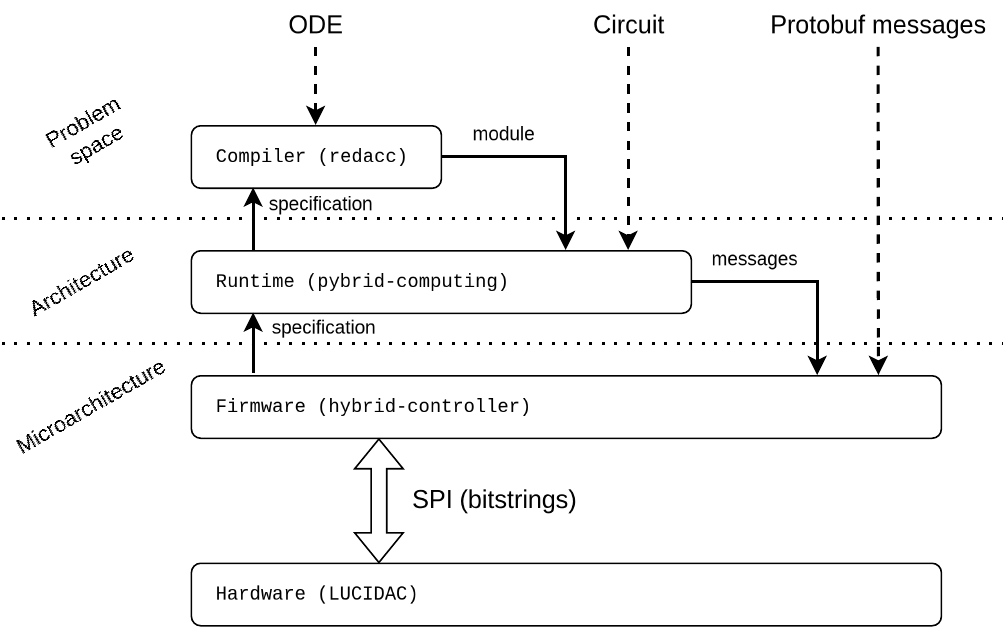}
    \caption{Overview of the LUCIstack across the problem-space, architecture, and microarchitecture layers. The compiler translates ODE specifications into circuit modules, the runtime manages circuit execution and protobuf-based commands, and the firmware maps these commands to the LUCIDAC hardware through SPI bitstrings. Upward-facing interfaces expose specifications between layers, while downward-facing paths carry compiled modules and runtime messages toward execution.}
    \label{fig:software_stack}
\end{figure*}
The availability of reconfigurable hybrid-analog digital computers, such as the REDAC platform, alongside emerging prospects for highly integrated analog computers on chips, requires a corresponding software infrastructure. Programming stored-program digital systems has become routine. Scripting languages, compilers, and debugging tools form an ecosystem that abstracts away hardware complexity. The wide adoption of digital computing partly rests on successive layers of abstractions that enable practitioners to solve problems without mastering the underlying machine architecture. For analog computers, no such layered approach exists. Users must reason directly about circuit structure and component behavior. Furthermore, modern computational workflows rarely consist of isolated single-step computations. Complex problems typically integrate multiple computational stages, each potentially optimized for different architectures. This demands tight integration of analog computing into established software stacks rather than treating it as an isolated accelerator.

In this section, we present our approach to addressing these software-related issues for both the REDAC and LUCIDAC platforms. The specifics reflect these platforms' hardware characteristics, though the underlying principles generalize. Parts of the software stack (LUCIstack) are open source, with comprehensive documentation available at \url{https://anabrid.github.io/lucistack/}.

The software challenges we address are not unique to analog computing but are present throughout the broader field of unconventional computing. Photonic and quantum computing, despite rapid advances in hardware design, face similar difficulties in establishing principled abstractions for their devices. In quantum computing, fundamental questions persist regarding whether qubits and gates constitute the optimal abstraction layer, or whether pulse-level control offers superior expressivity. Both fields must reconcile specialized hardware characteristics with the need for programmer-friendly interfaces.

Photonic processing currently targets specific application domains, particularly artificial intelligence acceleration and optimization via Ising machines.\cite{ahmed2025universal, hua2025integrated} In quantum computing, comprehensive software stacks like the Munich quantum software stack (MQSS)\cite{burgholzer2026mqss} and integration into high-performance computing facilities such as Oak Ridge National Laboratory\cite{shehata2025building} demonstrate that quantum computing has reached the systems integration phase. Reference~\onlinecite{burgholzer2026mqss} describes the MQSS stack using terminology familiar to computer scientists: architecture, microarchitecture, compiler frontend, intermediate representation. This reflects a common design principle across unconventional computing, namely to adhere to established abstractions from the stored-program paradigm while respecting the physical nature of the underlying hardware. Our approach follows this precedent, building on decades of experience in digital computing infrastructure.

%%%%%%%%%%%%%%%%%%%%%%%%%%%%%%%%%%%%%%%%%%%%%%%%%%%%%%%%%%%%%%%%%%%%%%%%%%%%%%%%%%%%%%%%%%%%%%%%%%%%%%%%%%%%%%%%%%
\subsection{A software stack for hybrid analog-digital computing}
%%%%%%%%%%%%%%%%%%%%%%%%%%%%%%%%%%%%%%%%%%%%%%%%%%%%%%%%%%%%%%%%%%%%%%%%%%%%%%%%%%%%%%%%%%%%%%%%%%%%%%%%%%%%%%%%%%
%
%
%
The software stack targets analog electrical computing while maintaining a general design applicable to other computational paradigms presented in this perspective. We pursue three goals derived from the requirements above. First, enable black-box operation where the device functions as a mathematical solver for directly representable problems (G1). Second, provide a consistent, stable, and flexible interface supporting low-level optimizations without sacrificing simplicity (G2). Third, integrate seamlessly into established workflows, conforming to user expectations from the digital, stored-program world (G3). We proceed from the top down, first defining mathematical abstractions for problems amenable to REDAC (\ie, the problem space), then describing the hardware structure representation (\ie, the architecture), and finally addressing physical constraints (\ie, the microarchitecture), as summarized in Fig.~\ref{fig:software_stack}.

%%%%%%%%%%%%%%%%%%%%%%%%%%%%%%%%%%%%%%%%%%%%%%%%%%%%%%%%%%%%%%%%%%%%%%%%%%%%%%%%%%%%%%%%%%%%%%%%%%%%%%%%%%%%%%%%%%

\subsubsection{Problem space}

The problem-type abstraction fundamentally shapes all layers of a software stack. In quantum computing, for example, choosing quantum circuits as the abstraction rules out the more general pulse-based approach used in analog quantum computing.\cite{meller2025analogqc} For analog computers, ODEs and PDEs have traditionally been used as the abstractions of choice. This view, however, overlooks devices such as the REDAC, which allow users to create feedback loops manually by connecting a multiplier's output to one of its inputs.

To include these types of problems, we define the class of problems $\mathcal{A}$ representable on REDAC. Let $F = \{f_a(t), g_b(t), \ldots\} \subseteq \mathbb{R}_{\ge 0} \rightarrow \mathbb{U}$ be a set of signals where $a, b, \ldots \in \mathbb{U}^n\,, n \ge 0$ are vectors of fixed parameters, and $\mathbb{U} = [-1, 1]$ corresponds to normalized voltages in $[-V_{\text{ref}}, V_{\text{ref}}]$ with $V_{\text{ref}} = 2\,\text{V}$ on the REDAC platform. Constraints $i \in I$ take the form
\begin{equation}
    \lim_{\Delta t \rightarrow 0} \frac{f(t + \Delta t) - \alpha_i f(t)}{(\Delta t)^{\alpha_i}} = F_i (t, f_a(t), g_b(t), \dots)\,,    
    \label{eq:redac_constraints}
\end{equation}
where $f \in F$ and $\alpha_i \in \{0, 1\}$. For $\alpha_i = 1$, Eq.~\eqref{eq:redac_constraints} reduces to ODEs. For $\alpha_i = 0$, feedback loops are modeled by exploiting the latency $\tau > 0$ of signals in real analog computers (typically nanoseconds on REDAC), creating a constraint $f(t) = f(t-\tau)$ that approximates a static algebraic relation. Any function $f \in F$ not represented by a constraint $i \in I$ is underdefined and must be an external input signal. By restricting signal values to $\mathbb{U}$, the set $\mathcal{A}$ contains only problems representable on REDAC-like computers. Exceeding $\mathbb{U}$ in practice results in ``overload'' mode, the handling of which is hardware-specific.

Each instance $a \in \mathcal{A}$ is designed to be solved with a single integration. In the following sections, we discuss solving problems from $\mathcal{A}$ individually whereas in applications, sequences of problems, where each instance's output influences the following instances, are required. While derived from the REDAC, this approach could be useful for the GPAC concept (see Sec.~\ref{sec:new_quest}), extending Shannon's original formulation by explicitly modeling feedback via physical latency rather than requiring strict polynomial ODE structure. Note that for practical modeling of ODE-like problems, the signal latency is no concern to users.
\subsubsection{Architecture}
The architecture layer represents hardware structure in an abstracted way, allowing users to create configurations that mirror actual problems from $\mathcal{A}$. This is analogous to the assembly level in digital computers, which provides low-level primitives through mnemonics such as \texttt{mov \%eax, 0x001}, while abstracting away details such as the physical register set, instruction pipelining, and virtual memory. For REDAC, this means building a structure that bundles firmware-abstracted components into semantic modules called ``blocks.''

The smallest building block of the REDAC is structured as follows: 
At its heart are 16 computing elements, typically eight integrators, four multipliers, and four constant value units, which are organized in M-blocks. Each of these elements expects input values as currents and yields output values as voltages. The sixteen output voltages are fed to the rows of a $16\times 32$ voltage-coupled switch matrix, the U-block, which acts as an expander. Since the values are represented by voltages, one output signal can be duplicated to multiple rows. 

Each of the $32$ output rows is connected to a digitally controlled coefficient which also performs an implicit voltage-to-current conversion. The resulting 32 currents are then fed to a second $32\times 16$ switch matrix. Since values are now represented by columns, multiple values can be implicitly summed in this current-coupled matrix, the I-block, which acts as a concentrator. Using this technique, no explicit summing computing elements are required, thus significantly simplifying implementation and programming of this system. The sixteen output currents of the I-block are then fed to the computing elements, thus closing an overall loop over all inputs, outputs, and computing elements. 

A small-scale system such as the LUCIDAC consists of one such cluster, while larger systems, up to a full REDAC, contain multiple clusters with additional interconnect structures allowing to pass multiple signals from one cluster to another.

A computational problem is now represented as a circuit by deriving configurations for the clusters of a given system. Routing signals, setting coefficients, and selecting signals to measure via ADCs is done digitally under control of an attached digital computer. 

%%%%%%%%%%%%%%%%%%%%%%%%%%%%%%%%%%%%%%%%%%%%%%%%%%%%%%%%%%%%%%%%%%%%%%%%%%%%%%%%%%%%%%%%%%%%%%%%%%%%%%%%%%%%%%%%%

\subsubsection{Microarchitecture}

The microarchitecture layer represents an implementation layer that abstracts away physical details while closely conforming to the architecture. Although analog computers are typically described in terms of the mathematical operations they implement, realizable systems remain subject to physical limitations. Among other constraints, signal and coefficient values are bounded by the capabilities of their electrical components. Unlike digital systems that can load/stop instructions, all operators within a problem (\eg, addition and multiplication) must be placed on hardware simultaneously. As no two electrical components have identical characteristics, production tolerances may occur and affect computation. There is also a small but non-infinitesimal signal delay on the chip. Finally, using voltages versus currents leads to different trade-offs as currents sum via Kirchhoff's law while voltages support fan-out natively, requiring conversion at appropriate points.

On this layer, firmware manages physical features through three abstractions. First, a hardware abstraction layer turns partial circuits into software-reconfigurable components with clean interfaces (\eg, multiple switching matrices into a single large matrix). Second, a calibration process measures internal gain and offset errors, compensating with software and hardware correction. Third, a control loop communicates with the host, receiving and applying circuit configurations while controlling integration and data acquisition. Effectively, the firmware transforms the analog computer into a data structure (\ie, abstracted components) and a collection of services. Technically, device-internal communications use the serial peripheral interface (SPI), while external host communication uses a binary protocol (see \url{https://github.com/anabrid/analog-protobuf} implemented in Google's protobuf\cite{currier2022protocol}).

The problem-space, architecture, and microarchitecture layers form a consistent paradigm, giving rise to a simple mechanism for mapping problems to implementations. By defining mappings between adjacent layers, these mappings can be composed to implement problems from $\mathcal{A}$ on REDAC machines.

%%%%%%%%%%%%%%%%%%%%%%%%%%%%%%%%%%%%%%%%%%%%%%%%%%%%%%%%%%%%%%%%%%%%%%%%%%%%%%%%%%%%%%%%%%%%%%%%%%%%%%%%%%%%%%%%%%
\subsection{Systems integration and optimization}
%%%%%%%%%%%%%%%%%%%%%%%%%%%%%%%%%%%%%%%%%%%%%%%%%%%%%%%%%%%%%%%%%%%%%%%%%%%%%%%%%%%%%%%%%%%%%%%%%%%%%%%%%%%%%%%%%%
\begin{lstlisting}[style=pythonstyle, caption={Analang DSL input for harmonic oscillator.}, label={lst:harmonic_dsl}]
let diff[h, t, t](t) = -h(t);
let h(t: 0) = 0.420000;
probe h(t);
\end{lstlisting}

\begin{lstlisting}[style=pythonstyle, caption={Lowered analogIR for the harmonic oscillator example.}, label={lst:harmonic_ir}]
module {
  analog.module @harmonic {
    %cst = analog.constant 0.0
    %cst_0 = analog.constant 1.0
    %0 = analog.time
    %c0_i64 = arith.constant 0 : i64
    %1 = arith.sitofp %c0_i64 : i64 to f64
    %2 = analog.fp_to_signal %1 : f64
    %3 = analog.mul %6, %cst_0
    %4 = analog.itor %6, %cst_1
    %cst_1 = analog.constant 4.200000e-01
    %5 = analog.mul %7, %cst_0
    %6 = analog.itor %7, %cst
    %7 = analog.neg %4
    analog.probe 0 = %4 * 1.000000e+00 + 0.000000e+00, ""
    %cst_2 = analog.constant 2.000000e-02
    %8 = analog.mul %0, %cst_2
    %cst_3 = analog.constant 1.000000e+00
    %9 = analog.neg %cst_3
    %10 = analog.sum %8, %9
    analog.plot(%10, %4)
  }
}
\end{lstlisting}

Harnessing the capabilities of the REDAC and larger integrated analog computers requires automating this mapping process, much as compiled languages do in the digital domain. To this end, the LUCIstack provides facilities for \emph{compiling} circuits and a framework for \emph{operating} the analog device from a digital host, thereby enabling tight integration into established workflows. Together, these components transform the analog computer from an isolated computing device into an application-specific accelerator. In Fig.~\ref{fig:software_stack}, we show how the compiler, runtime, firmware, and hardware components of LUCIstack realize the problem-space, architecture, and microarchitecture layers.

\subsubsection{The \texttt{redacc} compiler}
Mapping a problem from $\mathcal{A}$ requires assigning signals to lanes, operations to compute elements, routing the signals, and setting parameters and coefficients. The \texttt{redacc} compiler follows the standard front-end, middle-end, back-end structure. The front-end parses problems from $\mathcal{A}$ in the analang domain-specific language (DSL) that closely mirrors the mathematical notation and recursively transfers mathematical expressions into an intermediate representation (analogIR). In Listing~\ref{lst:harmonic_dsl}, we show an Analang implementation of a harmonic oscillator defined by its second-order time derivative and an initial value.

The middle-end uses analogIR (aIR), a dialect of MLIR~\cite{lattner2020mlir} expressing the circuit as a directed, cyclic graph using simple mnemonics similar to assembly language. Signal names are removed at this stage. Unique signals are represented by indices or labels starting with a \texttt{\%}. Compute operations (distinct from physical compute elements at this stage) use short mnemonics like \texttt{itor} for integrator. At this stage, the capabilities of the backend are not defined, hence operations not implemented through compute elements (\eg, \texttt{sin}) are allowed, and coefficient magnitudes and compute element counts remain unlimited. Therefore, analogIR may be considered an implementation of a \emph{virtual} device, much like LLVM IR.\cite{lattner2004llvm}

Once expressed in aIR, compiler passes progressively ``lower'' the code until the circuit structure corresponds one-to-one with the target device. For instance, each coefficient outside the range $[-8.0, 8.0]$ is decomposed by fanning out its associated signal across multiple lanes, distributing the coefficient among those lanes, and summing the resulting partial signals. Each pass produces valid aIR. In Listing~\ref{lst:harmonic_ir}, we show the final IR for the harmonic oscillator example of Listing~\ref{lst:harmonic_dsl}. Optional heuristic optimization passes leverage per-integrator configurable $k_0$ (ranging from $10^2$ on LUCIDAC to $10^4$ on REDAC, see Sec.~\ref{sec:intro}) to improve numerical accuracy within the limited hardware coefficient range.

Lastly, the compiler backend performs placement and routing. To this end, it decomposes the analogIR circuit graph into patterns that map onto the REDAC cluster architecture. Each pattern leads from a compute element or constant to another compute element, traversing the U-, C- and I-blocks. REDAC's architecture allows such pattern matches to map one-to-one onto paths inside clusters or between clusters through switching networks. The mapping routing is thus a constraint satisfaction problem (CSP) assigning trace endpoints to clusters, using Google's OR-Tools as solver.\cite{didier2023ortools} At this stage, a cluster is represented only by its capacity for compute operations, \eg, 8 integrators and 4 multipliers\textemdash no direct mapping between compute operation in the IR and physical compute elements exists.

After routing, actual compute elements are placed in clusters by mapping them to specific components. For each trace type, constraints are constructed, creating a CSP covering the intra-cluster mapping problem for all clusters simultaneously. The REDAC switch matrix for inter-cluster routing is limited. Every cross-cluster connection through a shared trace implies constraints for mapping traces onto physical lanes. Thus, a parallel generation of the mapping is only possible for connected components of the REDAC topology.

As a furter benefit, the compiler also addresses the inherent scaling problem for analog computers. An iterative loop compiles a problem from $\mathcal{A}$ into a circuit, executes it on the device, checks for an ``overload'' indicator (triggered when the $\mathbb{U}$ range is exceeded), then changes the function scale. This automatic scaling approach remedies precision limitations discussed in Sec.~\ref{sec:pde}.
\subsubsection{Operation}
\label{sec:operation}
The compilation process produces a REDAC device configuration consumable by the firmware. REDAC devices connect over ethernet, transmitting protobuf-based messages containing commands like \texttt{StartRun} or configuration and sample data. Any system interfacing with protobuf can integrate REDAC devices, allowing analog computing to be employed across various frameworks and programming languages.

The pipeline for implementing problems from $\mathcal{A}$ on an analog device enables multiple interaction methods at different granularity levels. Users can send circuits via command line interface and receive samples in bulk, which works well for script-based access. The entity object model provides a Python-based software twin of the device. Every configuration and operation on this model automatically issues to the actual device. Commands of the entity object model chain in a single session, with execution handed to the runtime for background processing of complex sequences including reconfigurations, integrations, and circuit evaluations. For maximum flexibility, custom code implementing the protobuf protocol grants full access at the architectural layer.

In summary, the presented LUCIstack provides abstractions at the problem-space, architecture, and microarchitecture levels, exposing the first two to developers. These abstractions, together with features such as calibration within the firmware,
implements black-box solver characteristics (G1). By using standards such as protobuf and Python-based integrations, LUCIstack
plugs in seamlessly into established digital environments and workflows (G3) while giving users fine-grained access to the device down to the architectural level. Lastly, besides the protobuf protocol, the session and entity object model concepts offer the
stable and flexible interface (G2) that can stay unchanged even when offloading session execution to on-device
digital processors. This abstraction stack enables direct implementation of problems from Sec.~\ref{sec:applications}. For example, the Kuramoto model's sine
coupling in Eq.~\eqref{eq:kuramoto} is handled by the nonlinear function generator, while the FPUT cubic terms in Eq.~\eqref{fput} require
coefficient decomposition across multiple lanes due to the $[-8, 8]$ constraint. For delay differential equations, one can leverage the hardware delay lines introduced in Sec.~\ref{sec:general_purpose}.

%%%%%%%%%%%%%%%%%%%%%%%%%%%%%%%%%%%%%%%%%%%%%%%%%%%%%%%%%%%%%%%%%%%%%%%%%%%%%%%%%%%%%%%%%%%%%%%%%%%%%%%%%%%%%%%%%%
\subsection{Perspectives}
%%%%%%%%%%%%%%%%%%%%%%%%%%%%%%%%%%%%%%%%%%%%%%%%%%%%%%%%%%%%%%%%%%%%%%%%%%%%%%%%%%%%%%%%%%%%%%%%%%%%%%%%%%%%%%%%%%

The initial version of LUCIstack is limited to a simple, synchronous circuit execution model with synchronous message exchanges between the controlling host machine and the REDAC device. As a first simplification, the stack implements a ``session protocol'' through which users can submit sequences of protobuf commands in advance for ordered execution, enabling large experiments. However, each exchange with the device creates latencies, which especially hurts batched execution of multiple circuits. Even small reconfigurations through partial configurations require a full host-device communication step.

The next major step for the stack would cut these latencies by enabling users to move session code to the digital part of the device while maintaining the session concept (G2). Program execution graphs would be sent to the device as custom C or C++ code that uses the session application programming interface directly on the device. Cutting these latencies enables advanced workflows beyond today's single-circuit model. Circuit-level decisions based on sampled values could trigger alternative configurations. This capability would be useful for mixed-integer optimization problems, in which the continuous component is mapped to analog compute elements and the decision variables are implemented digitally. Circuit optimization could then use generic search methods to identify circuits that represent a given problem based on data, akin to machine-learning approaches. Fast automated scaling would begin with an unscaled circuit, halt integration upon detecting an overload, identify the element causing the overload, and rescale the relevant coefficients. This capability would represent a significant advance in analog and hybrid-analog computing, broadening the range of applicable problems and enabling tighter integration with established workflows.
\section{\label{sec:applications}Applications}
Analog computation in hardware is useful for a wide range of research problems and practical applications. Analog computation realizes physical implementations of dynamical systems, so it is especially well-suited for modeling systems that are described by ODEs, PDEs, and stochastic differential equations (SDEs).\cite{killat2023solving,ulmann2026analog} One can also employ analog computation in optimization,\cite{hopfield1985neural,hopfield1986computing,jiang2023efficient,Kalinin2025,khoja2026automated} sensing and edge computing,\cite{1528676,yass} and circuit-native computing that is based on networks of minimalist nonlinear and chaotic oscillators.\cite{tanaka2019} 
There are also analog implementations of energy-efficient large language models\cite{leroux2025analog}, spiking neuronal-network models,\cite{wu2025dynamics}, and a variety of other systems.
\subsection{\label{sec:dynamical_systems}Dynamical-systems modeling---By Mason A.\ Porter, Marcelo J.\ Rozenberg, and Lucas B\"ottcher}
There has been a lot of research on applications of analog computers in dynamical-systems modeling.\cite{ulmann2026analog} In such implementations, one uses the continuous-time operation of analog computers to avoid the numerical time-stepping that is required by conventional digital solvers.\cite{koppel2021} 
In this subsection, we consider a variety of dynamical systems (\eg, several types of coupled-oscillator systems) that are not easy to implement and scale
on analog computers.
In other parts of this perspective article, we consider dynamical systems and differential equations that are easy to implement and scale on analog computers.
In Sec.~\ref{sec:general_purpose}, we discussed the wave equation, Burgers equation, and the canonical Lorenz system. In Sec.~\ref{sec:pde}, we consider the heat equation.

It is instructive to discuss some dynamical systems that are not straightforward to implement on analog computers. In particular, we consider (i) delay differential equations, (ii) coupled phase oscillators, (iii) lattices of nonlinearly coupled oscillators, and (iv) neuronal systems. It would be very interesting to implement them on analog computers, and we thus discuss some possible remedies that can hopefully facilitate it.

We start with delay differential equations. One can write a delay differential equation for $x(t)\in\mathbb{R}^n$ as the vector field
\begin{equation} \label{vf}
    \dot{x} = f(t,x,x_t)\,,
\end{equation}
where $x\equiv x(t)$ and $x_t(s) = x(t - s)$ with $s \in (0,\tau]$ denotes the recent history of a trajectory. The vector field~\eqref{vf} depends not only on the current time $t$ and state $x(t)$ of a system but also on its past states.
One initializes a standard integrator circuit with finite-dimensional initial conditions, whereas delay differential equations require an initial history function. The analog implementation of a delay differential equation thus requires additional mechanisms, such as delay lines (see Secs.~\ref{sec:gpac_building_blocks} and \ref{sec:operation}) or hybrid memory elements.\cite{losson1993solution}

Another important challenge is the analog simulation of systems of coupled oscillators, such as coupled phase oscillators. The most paradigmatic model of coupled phase oscillators is the Kuramoto 
model,\cite{kuramoto1975self,rodrigues2016kuramoto} in which the the phase $\theta_i(t)$ of oscillator $i\in\{1,\dots,n\}$ at time $t$ evolves according to
\begin{equation}
    \dot{\theta}_i = \omega_i + \frac{1}{n}\sum_{j = 1}^n K_{ij} \sin(\theta_j - \theta_i)\,,
    \label{eq:kuramoto}
\end{equation}
where $\omega_i$ is the intrinsic frequency of oscillator $i$ and $K_{ij}$ is the coupling strength between oscillators $i$ and $j$. In addition to the sinusoidal coupling term and the potentially many pairwise interactions, implementing \eqref{eq:kuramoto} in an analog circuit requires computing phase variables modulo $2\pi$ (or modulo some other constant). One possible remedy is to represent phases using complex variables 
$X_i(t) = \exp(\mathrm{i}\theta_i(t))$ or equivalently by using both their sine and cosine components. This avoids explicit phase wrapping, but it doubles the number of real state variables. An analog implementation of coupled Kuramoto oscillators was demonstrated recently in~Ref.~\onlinecite{freitas2026}.

A related interesting test case for analog computers are nonlinear lattice systems such as Fermi--Pasta--Ulam--Tsingou (FPUT) lattices.\cite{Fermi1955} The FPUT model takes the form of a large set of coupled ODEs.
In one dimension, the FPUT system consists of an array of harmonic oscillators that are coupled to each other by nonlinear springs. One of the original coupling functions is cubic and yields the FPUT-$\beta$ model
\begin{align}\label{fput}
	\ddot{x}_1 &= (x_2 - 2x_1) + \beta\left[(x_2 - x_1)^3 - x_1^3\right]\,, \notag \\
	\ddot{x}_i &= (x_{i + 1} + x_{i -1 } - 2x_i) + \beta\left[(x_{i + 1} - x_i)^3 - (x_i - x_{i - 1})^3\right]\,, \notag \\	
	\ddot{x}_n &= (x_{n - 1} - 2x_n) + \beta\left[- x_n^3 - (x_n - x_{n - 1})^3\right]\,,
\end{align}
where $x_i\equiv x_i(t)$ (with $i \in \{2, 3, \ldots, n - 1\}$), the lattice has $n$ oscillators, and $\beta$ is a coupling coefficient.

In the original FPUT numerical experiments,\cite{Fermi1955} Fermi et al.\ excited the dynamical system~\eqref{fput} in the first (\ie, lowest-frequency) mode of its linearization. They expected that the energy of the first mode would spread to higher modes and eventually lead to equipartition of energy among the modes.  However, they observed instead that the energy eventually returned to the first mode after initially spreading to higher modes. Based on these pioneering numerical computations, it seemed that such energy ``recurrence'' may continue indefinitely. Many researchers subsequently examined such recurrences,\cite{FPUT-review} and researchers have even observed very long ``superrecurrences''\cite{tuck-menzel} in numerical studies of FPUT systems.

A notable simplification in~\eqref{fput} is that all of the units are identical, but obviously that is never the case in real systems. In real life, the components of a system are not homogeneous, even when one attempts to create components that have the same properties. Instead, they are only ``the same'' within some tolerance, and such mismatches have fundamental effects on properties like energy recurrence.\cite{nelson2018,li2023} 
A key question for analog-computing systems is to implement FPUT recurrences and superrecurrences.

Models of neuronal systems\cite{ermentrout2010} (and of excitable systems more generally) will also benefit from analog implementations.
Even simple models of spiking neurons combine slow continuous-time evolution with fast threshold-triggered events and resetting, which yields stiffness in numerical computations. It has been noted that volatile memristors can describe voltage-gated ionic channels\cite{chua2012hodgkin} and provide useful circuit elements for neuromorphic hardware.\cite{delvalle2018challenges} There are recent analog implementations of spiking-neuron models that employ widely available and affordable electronic components.\cite{rozenberg2019,wu2023bursting} These implementations use a thyristor-based memristive circuit to reproduce neuronal dynamics, and there are also implementations 
of simple analog circuits with both excitatory and inhibitory biomimetic synaptic coupling.\cite{wu2025dynamics,Depannemaecker2026} 
These neuron and synaptic units provide building blocks for spiking neuronal networks with arbitrary connectivity and the ability to learn through synaptic plasticity.\cite{Hollande2026}
\subsection{\label{sec:pde}Partial differential equations (PDEs)---By Frederic Bergner}
\begin{figure}
    \centering
    \includegraphics[width=\linewidth]{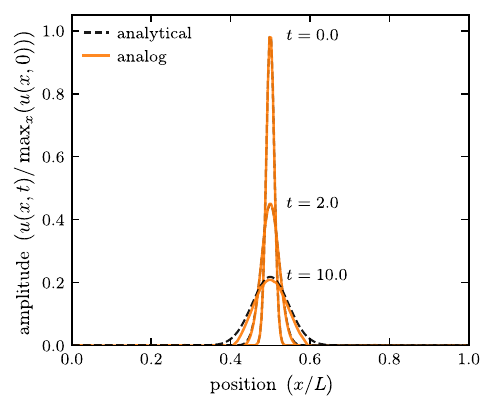}
    \caption{Solution of the one-dimensional heat equation. Diffusion of an initially Gaussian peak computed on the REDAC platform using 250 spatial nodes and continuous-time integration. Solid orange and dashed black lines indicate the REDAC (analog) and analytical solutions, respectively.}
    \label{fig:heat_equation}
\end{figure}

PDEs arise in many practically relevant fields, including fluid dynamics and finance.\cite{adler2025NumericalPartialDifferential} On digital computers, PDEs are typically solved by discretizing the spatial domain and numerical timestepping to integrate the resulting system in time. Time stepping may use either explicit or implicit methods. Explicit methods often require small time steps to resolve fast dynamics and maintain numerical stability. Implicit methods can remain stable at larger time steps, making them suitable for capturing slower dynamics in stiff systems, but each step is computationally more expensive. Analog computers perform time integration directly in the continuous-time domain. Continuous-time integration resolves both fast dynamics\textemdash without stability constraints of explicit time stepping\textemdash and slow dynamics\textemdash without the computational overhead of implicit methods. This offers significant advantages for PDEs that exhibit a wide range of temporal scales and lead to stiff ODE systems. An important representative of this type of PDE are the Navier--Stokes equations. Especially for technical applications, the spatial and temporal scales can span several orders of magnitude. Therefore, time stepping becomes expensive, forcing a trade-off between resolving fast dynamics and neglecting or modeling them. Because continuous-time integration circumvents these problems entirely, analog computing offers benefit for this class of problems.
Spatially discrete, continuous-time analog formulations of PDEs have been applied in areas such as acoustics and fluid dynamics.\cite{malavipathirana2025AnalogComputingNonlinear, mandal2024IntegratedAnalogComputers} The continuous-time solution of the one-dimensional heat equation, a PDE that arises frequently in many areas of physics, is shown in Fig.~\ref{fig:heat_equation}. The continuous-time integration of the diffusion of a Gaussian peak is performed on the REDAC platform (see Sec.~\ref{sec:redac_platform}) on 250 spatial nodes and 
compared with the known analytical solution. Related analog approaches have also been used to generate stochastic trajectories required by PDE-based models in finance, such as the Black--Scholes--Merton framework.\cite{huang2018CaseStudyAnalog} More broadly, prototype analog and hybrid analog-digital chips have demonstrated reduced solution times and energy consumption relative to digital microcontrollers for selected low-precision differential-equation workloads.\cite{guo2016EnergyEfficientHybridAnalog}

Despite the advantages of continuous-time integration, analog computing is not widely used in high-performance computing because of its lower precision relative to digital computing, susceptibility to thermal noise and hardware imperfections, and hardware requirements that often grow with problem size. For the comparatively simple case of the one-dimensional diffusion of a Gaussian peak, depicted in Fig~\ref{fig:heat_equation}, the deviation from an analytical solution increases with runtime as errors due to hardware imperfections and thermal noise accumulate over time. Achieving higher precision generally requires disproportionate increases in energy consumption and silicon area. These costs can grow exponentially with the number of effective bits. Analog circuits may therefore offer advantages in energy consumption and silicon footprint for low-precision calculations, whereas digital circuits generally scale more favorably for high-precision computation.\cite{sarpeshkar1998AnalogDigitalExtrapolating} Amplified thermal white noise can, on the other hand, be leveraged to solve stochastic PDEs by physically introducing stochastic forcing into the analog circuit.\cite{huang2018CaseStudyAnalog} While the spatially discrete, continuous-time approach avoids certain problems associated with time discretization, spatial discretization results in large systems of coupled ODEs. The time integration of each ODE requires an analog integrator, along with summers and multipliers. This one-to-one scaling poses challenges for real-world applications, as the resulting systems may contain millions of ODEs and require analog computers of comparable scale. One way to scale analog circuits is through miniaturization.\cite{kinget2015ScalingAnalogCircuits} Because analog components dissipate relatively little energy, advanced manufacturing techniques may enable three-dimensional analog circuits. Such architectures are being developed primarily in the context of neuromorphic computing\cite{adam20173DMemristorCrossbars} and could accommodate more analog components on a single chip while also enabling more efficient interconnections.
\subsection{Optimization---By Lucas B\"ottcher}
Solving an optimization problem (\ie, maximizing or minimizing an objective function, often subject to constraints) is closely related to finding the ground state of a physical system.\cite{kirkpatrick1983optimization} In many physics-inspired approaches, the objective and constraints are encoded in an energy function, so that low-energy states correspond to good candidate solutions.

The 2024 Nobel Prize in Physics was awarded to John J.\ Hopfield and Geoffrey Hinton ``for foundational discoveries and inventions that enable machine learning with artificial neural networks''. Beyond the connection between artificial neural networks and magnetic systems such as the Ising model,\cite{bottcher2021computational,bottcher2025statistical} the scientific background document of the Nobel Committee for Physics\cite{nobel2024physics_background} also highlighted the close correspondence between mean-field approximations of binary systems and continuous-time dynamics implemented in analog computers to search for low-energy configurations that encode candidate solutions to combinatorial optimization problems.
\begin{formal}
Optimization occurs through integration of the equations of motion of an electronic circuit, during which the nodes evolve without instructions from a central unit. This approach constitutes a pioneering example of using a dynamical system to seek solutions to difficult discrete optimization problems.\cite{mohseni2022ising}
\vspace{1em}\\
\footnotesize{The Nobel Committee for Physics (Scientific Background to the Nobel Prize in Physics 2024)}
\end{formal}

This line of work goes back to Hopfield and Tank,\cite{hopfield1985neural,hopfield1986computing} who proposed networks of nonlinear analog neurons for combinatorial optimization. One of the examples that Hopfield and Tank considered was the traveling salesman problem. Subsequent neural mean-field methods extended related ideas to other problems, including knapsack problems.\cite{ohlsson1993neural,ohlsson1997study} A survey of this early generation of neural optimization methods is available in Ref.~\onlinecite{smith1999neural}.

Modern analog and hybrid optimization approaches include Ising machines,\cite{mohseni2022ising} which aim to find exact or approximate ground states of Ising-type energy functions. Oscillator-based Ising machines are closely related to the Kuramoto model~\eqref{eq:kuramoto}, where binary variables can be represented by oscillator phases, such as phases 0 and $\pi$. This phase-based representation has historical roots in the parametron work of Goto\cite{4065825} and related work by von Neumann\cite{von1957non,4065705}.
\subsection{\label{sec:edge_computing}Edge computing---By Bhaskar Choubey}
Analog computing at the sensor edge has a long history of research, dating back to the earliest developments of transistors. Since most physical stimuli are inherently analogue, converting them into an electrical signal typically requires an analog circuit. Even in the simplest sensing systems, analog computation is provided in form of buffering and signal conditioning before digitization. Each stage of this sensing pipeline\textemdash transduction, conditioning, and digitization\textemdash has been investigated to perform early stage analog computation to improve the efficiency and accuracy of downstream digital algorithms while simultaneously reducing communication bandwidth and energy consumption. Such analog computation is particularly attractive at the sensor edge as the edge is typically constrained by stringent power, area and bandwidth budgets. Since data converters and transceivers are among the most energy-intensive components of many sensing systems, reducing or postponing digitization can provide substantial improvements in overall system efficiency. At the same time, proliferation of low-cost sensors has led to exponential rise in data being produced by these sensors. This is hence adding new bottleneck of sensor-processor and sensor-memory in addition to memory-processor bottleneck in the traditional pipeline of sensory data processing. Once again early analog processing at the sensor site could alleviate some of these problems. Existing approaches to analog computation at the sensory edge can be broadly classified into three categories according to where and how the computation is performed. 
\subsubsection{Computing in sensors}
The first category, ``computing in sensors'', integrates computational functionality directly within the sensing device. Early pioneering work by Mead introduced ``neuromorphic sensors'',\cite{58356} particularly logarithmic retinal image sensors capable of adapting to changing illumination levels. It inspired an extensive body of research into neuromorphic vision sensors that emulate the response characteristics of the human visual system.\cite{222815,6861463,choubey_aic} Unlike conventional image sensors with linear intensity responses, biological vision exhibits approximately logarithmic or Stevens' power-law behavior, enabling operation over an exceptionally wide dynamic range.\cite{leeAdaptive2025,4268351,Collins2008} Circuits using the subthreshold operation of metal oxide semiconductor (MOS) transistors have been proposed to realize compact nonlinear logarithmic response circuits. Beyond simple logarithmic compression, circuit designs have also been proposed which implement complex functions like tone mapping operators,\cite{Gouveia} essential for displaying high-dynamic-range images on conventional displays with limited dynamic range. Although these operations are computationally demanding in the digital domain, they can often be approximated efficiently using compact analog circuits embedded within individual pixels. The same design philosophy of computation using signal reconditioning has also been utilized beyond optical sensing including pH and chemical sensors\cite{1232939,9471804} as well as infrared detectors\cite{budzier,344700}. In addition to response shaping, early stage computation has also included calibration and noise reduction using techniques like correlated double sampling to suppress fixed-pattern noise by compensating for device mismatch.\cite{Gouveia1}
\subsubsection{Computing using sensory elements}
The second category, ``computing using sensory elements'', repurposes sensing devices themselves as computational primitives. Rather than acting solely as sensors, devices such as photodiodes, microelectromechanical system (MEMS) resonators, or specialized transistors have been re-configured to perform general-purpose analog computation. As an example, MOS-based chip fabrication processes provide optical diodes which have been used as the foundation of optical computing,\cite{9064516,ocm} while single photon avalanche devices have provided platforms to undertake quantum computing.\cite{11204475} Alternatively, MEMS resonators have been used as foundation block for different computing tasks like analog multiplexing,\cite{G1,3_1} reservoir computing\cite{9699777,11284744} as well as a variety of neural networks.\cite{8443098} Resonant MEMS devices typically solve the linear pencil second order differential equation and have been used in a range of sensory applications including gyroscopes, mass, velocity, pressure and temperature.\cite{beeby2004mems} Reducing their dimensions, while increases their sensitivity, has also known to introduce nonlinearity, which provided further opportunity for mechanical computing. In addition, resonant sensors have also been coupled for sensory enhancement, quality factor control as well as mechanical multiplexing. These arrangements inherently solve an array of coupled linear or nonlinear differential equations. In their simplest form, the collective behavior of these resonant sensors provides mechanical multiplexing.\cite{G2} However, they can also be envisioned as mechanical computers with potential application in harsh environments.
\subsubsection{Computing for sensors}
The third category, ``computing for sensors'', performs the initial stages of signal processing immediately adjacent to the sensing element. This approach has attracted considerable attention with the recent growth of edge artificial intelligence. Examples include analog implementations of convolutional filters, edge detection, feature extraction, object detection, and the early layers of neural networks integrated directly with sensors, often for autonomous driving.\cite{1528676,yass} Diverse sensory systems have been reported with early signal processing with some examples including gas sensing front-ends,\cite{7938702} auditory sensors including cochlear implants,\cite{102813,8048493} tactile sensors\cite{7784764} as well as olfaction\cite{4061019}. This rapidly increasing domain of computation, nevertheless, can be considered to be a subset of the wider analog computing being discussed in this perspective, implemented under constraints of limited power, area and speed requirements. 

\begin{figure*}
    \includegraphics[width=0.9\textwidth]{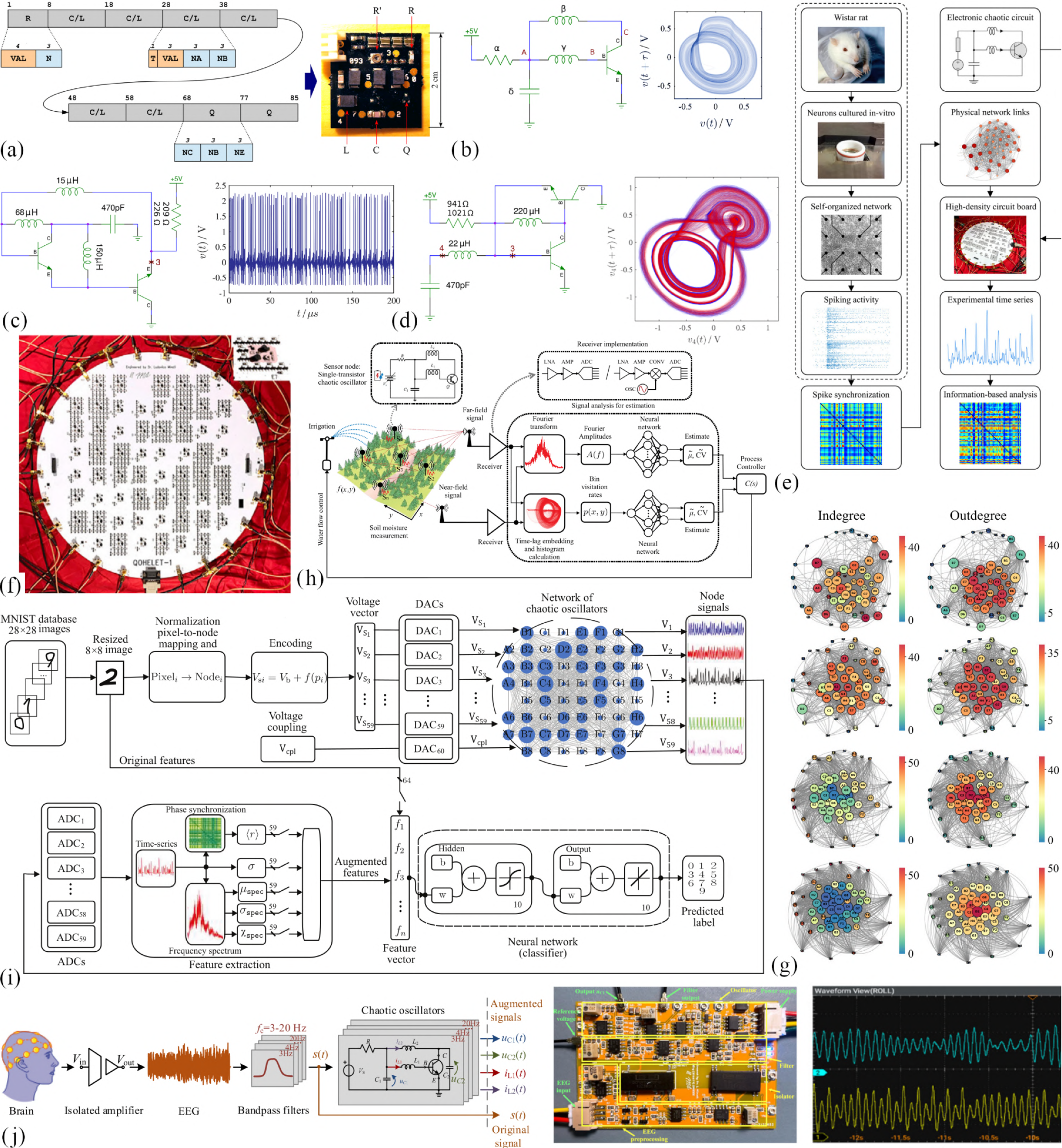}
    \caption{Overview of analog computing by means of elementary, transistor-based circuits. (a) Representation of circuit topology as bit-string, suitable for large-scale searching of new arrangements. (b) Simplest single-transistor configuration, (c) spiking configuration, and (d) double-scroll configuration. Panels (a)--(d) are reproduced from Ref.~\onlinecite{minati2017atypical} by permission of the American Institute of Physics.
     (e) Conceptual framework for the design of an extended physical reservoir network based on a biological neural culture. (f) Resulting printed circuit board realization. (g) Diversity of emergent directed information flows between hubs and the rest of the network. Panels (e)--(g) are reproduced from Ref.~\onlinecite{minati2025culture} by permission of Elsevier. (h) Hypothetical scenario of distributed sensing in an agricultural internet-of-things application, and its testing in a benchtop experiment. Reproduced from Ref.~\onlinecite{zhao2026distributed} by permission of Elsevier. (i) Application of the reservoir network to handwritten digit recognition, through mapping the intensities of the most representative image pixels to the control voltages of the most central oscillators. Reproduced from Ref.~\onlinecite{tagnefossi2026digits} by permission of Elsevier. (j) Analog-domain integration of an EEG amplifier and single-transistor chaotic oscillator, demonstrating the generation of partial synchronization and nontrivial responses to the input. Reproduced from Ref.~\onlinecite{fu2025eeg} by permission of Elsevier.
     }
    \label{fig:circuit_based}
\end{figure*}

\subsection{Computing using networks of minimalist chaotic circuits---By Ludovico Minati and Manyu Zhao}
In this section, the possibility of using elementary chaotic circuits to realize emergent forms of analog computing is considered in detail, and some representative results are surveyed.
\subsubsection{Circuit-native nonlinear dynamics as the substrate of physical computation}
Analog computer architectures are usually based on ODE models, in which integrators, summers, multipliers, and nonlinear function modules are interconnected to realize a predetermined vector field, with the mathematical formulation directly determining the hardware design. On the other hand, physical reservoir computing leverages existing real-world properties of materials, systems and collectives, in order to harvest their capability to do computation based on their inherent characteristics, seeking instead to transform the input and output information in order to unlock the computational capability towards a given application.\cite{tanaka2019} Since analog electronic circuits are physical systems, the same conceptual rule can also apply. In this framework, a noteworthy example is that of chaotic oscillators, which are particular arrangements, often discovered serendipitously, wherein the design is largely determined by ab-initio device and interconnection choices, and by the physical characteristics of the network (the network structure, the Kirchhoffian properties, and the behavior of the semiconductor devices). In other words, the differential equations are not given a priori but have to be derived a posteriori from the circuit structure, and reduced-order models are only derived subsequently. While a comprehensive review of this approach is beyond the scope of this work, a noteworthy case is that of atypical circuits based on one or at most two transistors, surrounded by a limited number of inductors and capacitors realizing reactive networks, and a resistor placed in series to a DC power supply. In the absence of a fully systematic design methodology, such topologies can be explored through extensive numerical searches, based on genetic algorithms or random trajectories (see Fig.~\ref{fig:circuit_based}a).\cite{minati2017atypical}

The simplest known structure obtained using this approach comprises one bipolar transistor, one resistor, two inductors, and one capacitor (see Fig.~\ref{fig:circuit_based}b). Following its discovery and experimental realization, the behavior of this circuit was described as a relaxation oscillator in a reciprocal interaction with a damped resonator, with chaotic dynamics arising due to the interaction between the two.\cite{minati2022mechanisms} The same design space was also found to contain a dual-transistor spiking oscillator, which generates spike trains resembling neural bursts while being devoid of any neuromorphic aspects (see Fig.~\ref{fig:circuit_based}c), and a six-component double-scroll generator, where the incomplete symmetry in connection between two transistors and alternation of saturation to cutoff leads to competing switching behavior (see Fig.~\ref{fig:circuit_based}d).\cite{minati2022spiking,cen2025doublescroll} Similar approaches could also be applied to other analog implementations that are better suited for integrated circuits. One notable example is a cell based on current-starved, cross-coupled CMOS ring oscillators of lengths 3, 5, and 7, interconnected through diodes. This architecture can generate continuous chaotic signals, as well as spiking, bursting, and quenched dynamics, over a frequency range spanning at least four orders of magnitude, while maintaining low area and power consumption.\cite{minati2019cmos}

In other words, in the equation-based approach, the designer is directly involved in choosing the dynamical system and then mapping it to hardware using operational amplifiers, multipliers, integrators and the link. In the circuit-native approach of physical computing, the designer chooses the device, typically an elementary element such as a transistor, the network topology, biasing conditions, coupling schemes, and observables, but the vector field is not automatically assigned. The resulting system computes by using the transformations and dynamical behavior that are inherent to the physical substrate, not by forcing it to follow a mathematical model. The first approach is evidently more general, but leads to circuits that are large (up to thousands of active junctions), power-consuming and slow, typically operating in the kilohertz region. The second approach is less general but leads to highly compact circuits that are readily scalable to the gigahertz range and exhibit surprisingly rich dynamics.

\subsubsection{Collective generativity and neural-like emergent behaviors}

The compact nature of transistor-based chaotic oscillators lends itself well to the realization of extended networks, which can be problematic in all but the largest conventional analog computers. The combination of these elementary circuit elements leads to a strikingly rich dynamical repertoire that exceeds what one might expect given the low dimensionality of the node dynamics. For example, a 59-node network with structural connectivity replicating a dissociated neuronal culture was recently found to exhibit emergent phenomena, including spontaneous dynamical differentiation, waveform complexity modulated by node degree, and asymmetric information flows, all arising despite the use of simple diffusive symmetric physical couplings (see Fig.~\ref{fig:circuit_based}e,f). In particular, depending on the coupling configuration, both hub and peripheral nodes may self-organize to act as sources or sinks of information (see Fig.~\ref{fig:circuit_based}g).\cite{minati2025culture}
On the same physical substrate, hierarchical cluster synchronization also emerged in the absence of deliberately imposed network symmetries, and this organization persisted reliably even in the presence of unavoidable device mismatches, non-idealities, and noise.\cite{manjunatha2026cluster} Posterior analyses using partial-information decomposition further revealed the coexistence of redundant and synergistic informational contributions, with synergy becoming particularly pronounced even among node pairs lacking direct structural connectivity.\cite{bara2026synergy} Experiments conducted systematic lesioning of the network according to trajectories following increasing or decreasing centrality, further demonstrating that the degradation of network function is strongly topology-dependent. Damage propagating from the periphery to the core resulted in gradual deterioration, while core-to-periphery lesions precipitated an abrupt collapse of both synchronization and graph-based organization, typically preceded by a peak in modularity.\cite{tegnitsap2026lesions} At a larger scale, a more realistic model of a low-resolution human connectome was designed in the form of a low-power integrated circuit based on the aforementioned CMOS circuit, recovering nontrivial correspondences between structural and functional connectivity, with network hubs exerting a pronounced influence on the collective dynamics.\cite{kuate2025connectome}

The noteworthy point here is that these phenomena align the phenomenology of brain dynamics closely as observed through neurophysiological recordings, but this does not imply that the circuits are neuromorphic: evidently, they are not, since they implement neither ionic channels nor a specified neuron model, and are in fact considerably smaller than typical realistic electronic mirrors of physiology. Instead, cluster formation, symmetry breaking, source-sink differentiation, critical transitions, and multivariate information structure recur because nonlinear dynamics possess universal organizational mechanisms that are concordant across materially unrelated systems. These networks are therefore physical experimental models of generic nonlinear organization, not transistor-level replicas of neural tissue, and as such, exemplify the profound generative potential of these circuits. As said, this potential is best expressed over extended networks of nodes, which may be more problematic to realize when each node requires several operational amplifiers, multipliers, and the like.\cite{minati2018silicon}
\subsubsection{Engineering relevance and applications}
Owing to their inconspicuous size and convenient scalability across frequency and power levels, transistor-based chaotic oscillators and networks have received increased attention for diverse engineering applications, where they can fulfill useful signal processing, augmentation, and integration functions without explicitly realizing particular equations.\cite{minati2024review} They represent, therefore, a particularly deployable form of analog computing, all way down to directly driving, for example, industrial processes when the circuit itself is built using power electronics, replacing a separate numerical controller.\cite{minati2025manganese}

For example, in an initial optical proof-of-concept experiment, photovoltaic cells power single-transistor oscillators and transduced luminous flux into a control parameter determining the transition to chaos, with mutual inductive coupling between adjacent nodes enabling the ensemble dynamics to encode spatially distributed inputs. This configuration permitted estimation of the global average without necessitating interrogation of each node.\cite{minati2020distributed} Subsequent numerical investigations extended this architecture to more realistic irregular two-dimensional networks, demonstrating that both the mean and variability of local sensor readings could be inferred from near- or far-field antennas attached to a spectrum analyzer, and further explored the use of frequency-division and time-division multiplexing as scalable coupling strategies.\cite{tagnefossi2025distributed} Later, experimental validation with a 59-node benchtop system using the same network mentioned above confirmed that physical estimates of both the mean and the coefficient of variation could be obtained from either the frequency spectra or the time-delay reconstructions of the collective output (see Fig.~\ref{fig:circuit_based}h).\cite{zhao2026distributed} On the whole, these results illustrate the extent to which sensor signal processing functions could be integrated within the physical layer itself, in a distributed manner.

Besides statistical estimation, the same physical substrate can also be used for classification tasks, and a recent paper demonstrated its application to handwritten digit recognition. By using topology-aware input assignments and operating at an intermediate coupling strength, the network of single-transistor oscillators was used to transform low-resolution input images into class-dependent spatiotemporal responses (see Fig.~\ref{fig:circuit_based}i).\cite{tagnefossi2026digits} On the other hand, for temporal data, it has been shown that human electroencephalographic (EEG) signals could be fed, in real time and in the analog domain using a coupled pre-amplifiers, directly to the single-transistor oscillator, obtaining a response is partially phase locked but partially autonomous, leading to the generation of time-frequency features that are influenced by the structure of the input in a highly non-trivial way (see Fig.~\ref{fig:circuit_based}j).\cite{fu2025eeg} When applied to patient EEG data, these physical features improved seizure classification appreciably, heralding potential future use in near-sensor and in-sensor data augmentation for biomedical applications, where extracting information at low power from complex signals is particularly important.\cite{tabekoueng2026seizure}
\subsection{\label{sec:ai_applications}Transformers and attention---By Paul Manea}
Contemporary artificial intelligence applications are predominantly based on transformers because their underlying self-attention mechanism is able to capture long-range dependencies very effectively.\cite{vaswani2023attentionneed,lin_survey_2022} In many applications, transformers are trained autoregressively. That is, they learn to predict future tokens sequentially. From a hardware perspective, however, autoregressive inference is associated with substantial memory-movement constraints, limiting the performance of modern GPUs such as the NVIDIA H100.\cite{luo2024benchmarkingdissectingnvidiahopper,nvidiaH100datasheet} At each decoding step, the query must access the growing key-value cache of all previous tokens. This means that data must be transferred repeatedly from high-capacity, high-bandwidth memory to on-chip memory and finally to the streaming multiprocessors for computation. To reduce the cost of repeatedly moving attention data between memory and compute units, one can employ the in-memory computing approaches described in Sec.~\ref{sec:memristive}.

Attention introduces two additional hardware challenges. First, the score vector $QK^T$ scales with context length, such that digitizing all intermediate scores creates an ADC bottleneck that can offset the area and energy advantages of analog compute-in-memory.\cite{wan2022compute} Second, standard softmax is difficult to implement efficiently in hardware because it combines element-wise exponentiation with a global normalization across all tokens. Recent hardware-software co-design approaches consequently either replace softmax with simpler nonlinearities or reformulate where normalization is performed.\cite{FinkbeinerSoftmax2026}

To address these constraints, Ref.~\onlinecite{leroux2025analog} proposes an analog in-memory computing approach for transformer attention based on gain-cell memories. The main idea is to reduce key-value-cache movement and avoid intermediate analog-to-digital conversion. This is achieved by storing the dynamic key-value cache directly inside the memory arrays and performing the two dominant attention dot products, $QK^T$ and $\phi(S)V$, in the analog domain. Queries are encoded as pulse-width-modulated signals. Intermediate results are then processed by charge-to-pulse circuits and a HardSigmoid nonlinearity is applied between the two stages. 

The energy consumption and latency of the proposed analog approach is reduced by up to four and two orders of magnitude, respectively, compared with GPUs.\cite{leroux2025analog} At the same time, the text-processing performance is comparable to GPT-2.
\section{Summary and analog supremacy---By Ludovico Minati}
\label{sec:summary}
The computational-complexity profile of a GPAC can make physical evolution time independent of system dimension, at the cost of additional hardware occupancy. For example, consider an initial-value problem of the form $\dot{x}=f(x,t)$, with $x\equiv x(t)\in\mathbb{R}^{n}$, over a fixed simulated interval $T$, and let $M(n)$ denote the number of arithmetic operations required for one evaluation of $f$. A digital solver has to calculate $N_{\mathrm{s}}$ discrete updates, which require at least $n$ state-variable reads and writes per update, implying a state-memory traffic of at least $\Omega(nN_{\mathrm{s}})$, so that evaluating the solution costs $\Theta(M(n)N_{\mathrm{s}})$ arithmetic operations. Generally, one has $M(n)=\Theta(n)$ for bounded-degree sparse or locally coupled systems and $M(n)=\Theta(n^{2})$ for dense pairwise coupling. An implicit solver may additionally require repeated solution of an $n\times n$ linear system, with a dense direct factorization costing $\Theta(n^{3})$ per time step, although some methods can alleviate this cost.\cite{hairer1996solving} In a GPAC realization, the $n$ state variables are assigned to $n$ physical integrators and the operations defining $f$ to $O(M(n))$ concurrently operating circuit elements. Provided that the necessary resources are available and that no limitations require a reduction of the integration-rate factor $k_{0}$, the physical operating time is $T/k_{0}$ and, crucially, does not increase with $n$. In other words, increasing the system dimension increases circuit size rather than the number of sequential computational steps. Sec.~\ref{sec:redac_platform} provides an example whereby enlarging the spatial discretization of the Burgers equation led to elevated hardware occupancy and, at the same time, a higher speed ratio relative to the digital processor.

Stiffness leads to a second scaling issue. For a dynamical system whose linearization, described by its Jacobian (stability) matrix, has eigenvalues that all lie on the negative real axis and have widely separated magnitudes, an explicit numerical method with time step $h$ must satisfy a stability criterion $h\lesssim C/|\lambda_{\max}|$, where $\lambda_{\max}$ and $\lambda_{\min}$ denote the eigenvalues of largest and smallest magnitude, and $C$ is a constant of order unity that depends on the method (for instance, $C=2$ for the explicit Euler scheme). In situations where the simulation length is determined by the span of the slowest fluctuations, such as $T=O(1/|\lambda_{\min}|)$, the required number of steps grows with the stiffness ratio $|\lambda_{\max}|/|\lambda_{\min}|$. The use of implicit methods can remove or substantially relax this constraint; however, such methods still require accuracy-controlled time steps and nonlinear iterations or Jacobian solves. By contrast, an analog integrator has no discrete numerical time step, and hence no finite step count, which is the basis of the advantage for stiff systems discussed in Sec.~\ref{sec:pde}. Successfully integrating systems with fast-slow dynamics, nevertheless, does not come for free in engineering terms, because realizing widely different time constants in analog computers requires adequate circuit bandwidth and dynamic range, together with suppression of parasitic issues such as capacitor non-idealities. All this may eventually force choosing small values of $k_{0}$ to retain accuracy. Thus, analog integration can deliver a physical integration time independent of $n$, and consequently avoid the stiffness penalty; however, its total complexity in terms of resources is independent neither of $n$ nor of the required timescale separation. This is a key consideration, and it is precisely why the continuous-time complexity theory reviewed in Sec.~\ref{sec:new_quest} measures computational effort by the length of the solution trajectory rather than by physical time: rescaling time leaves the length of the solution curve unchanged, so that trajectory length, unlike physical duration, provides a meaningful analog counterpart of digital time complexity.\cite{JournalACM2017}

Computational completeness and conditioning also need to be considered. The former pertains to which solutions can be attained through a given computational process. The latter is more specifically about how strongly errors in the specified problem are amplified in the process of reaching a solution. Sec.~\ref{sec:new_quest} established that, under the convergent-output definition, whereby an input is represented by computable data and the output approaches the requested value to a prescribed accuracy, GPAC-computable real functions coincide, on compact domains, with Type-2 Turing-computable real functions.\cite{JOC2007} In other words, analog computation is neither less nor more complete than digital computation. Models such as the well-known Blum--Shub--Smale machine, which operate on exact real numbers and count arithmetic operations at unit cost, are not at odds with this conclusion: they represent an abstract situation which is fundamentally not available to physical hardware.\cite{BlumShubSmale1989} Importantly, the equivalence in computability does not, however, imply equivalence in the resources required to obtain an accurate result. For a problem $y=F(d)$, where $F$ maps the input data $d$ to the solution $y$, the local relative condition number $\kappa_F(d)$ quantifies how strongly a small relative perturbation $\Delta d$ of the data is amplified into a change $\Delta y$ of the solution: to first order, $\frac{\|\Delta y\|}{\|y\|}\lesssim\kappa_F(d)\frac{\|\Delta d\|}{\|d\|}$. Accordingly, for a realization operating with an internal error on the order of $2^{-q}$, obtaining $p$ reliable output bits generally requires $q\gtrsim p+\log_2\max\{1,\kappa_F(d)\}$. This requirement is dictated by the problem itself, not by the technology, and it therefore applies whether the solver is implemented as a numerical algorithm or as a circuit in the analog domain.\cite{Rice1966Condition} In the former case, a target level of $q$ can be reached by making recourse to longer floating- or fixed-point formats, tighter solver tolerances (finer steps), and so forth. In the latter case, better component accuracy and quality, higher signal-to-noise ratio, tighter calibration, the use of averaging, and higher converter resolutions are required. While effective at reducing errors in practice, neither the bespoke hardware calibration and scaling techniques covered in Sections~\ref{sec:general_purpose} and~\ref{sec:software}, nor the use of various preconditioning methods, can entirely eliminate this sensitivity.\cite{feinberg2021analog,zuo2025precise} In summary, under reasonable robustness assumptions, analog and digital computers are equally complete, and no implementation on either side can escape the working precision dictated by the conditioning of the problem. On the analog side, this precision is a computational resource in its own right: as reviewed in Sec.~\ref{sec:new_quest}, it plays for continuous-time models the role that memory (space) plays for Turing machines.\cite{Icalp2024}

Taken together, Sections~\ref{sec:intro}, \ref{sec:new_quest}, and~\ref{sec:bio_computing} support three important conclusions regarding the nature of analog computation. First, Turing machines and GPACs profoundly differ in the way computation is represented and executed, but not in what can be computed. This is a point that needs to be underlined. Namely, a Turing machine steps through discrete states and eventually returns an answer if it reaches a designated halting state. On the other hand, a GPAC evolves along a continuous trajectory defined by polynomial ODEs, and its output is either a signal generated at physical time $t$ or a state variable converging to the requested value. But this difference does not influence what can be computed. Second, the halting limitation also applies to polynomial-ODE computation. Polynomial ODEs can simulate a universal Turing machine, even in the presence of bounded perturbations. It follows that a hypothetical device capable of deciding whether the trajectory of an ODE eventually reaches a target accepting region of the state space would be able to decide the Turing halting problem, which is impossible.\cite{Tur36,GracaCampagnoloBuescu2008} Third, the biological computation described in Sec.~\ref{sec:bio_computing} is unquestionably analog at the level of physical implementation because molecular concentrations, reaction rates, diffusion fields, stochastic events, and internal states are physical in nature. However, biology also differs from a GPAC because plasticity, adaptation, degradation, and evolution can continuously modify the computation being processed, something that conventional electronics cannot do well. In the future, memristive computing as described in Sec.~\ref{sec:memristive} and the use of other adaptable platforms may bridge this gap. It should, however, be noted that while biological computation is an adaptive physical analog process, it is not a more computationally complete one: continuous chemical reaction networks, for instance, are Turing-complete rather than super-Turing.\cite{CMSB17vantardise}

\begin{figure}
    \centering
    \includegraphics[width=\linewidth]{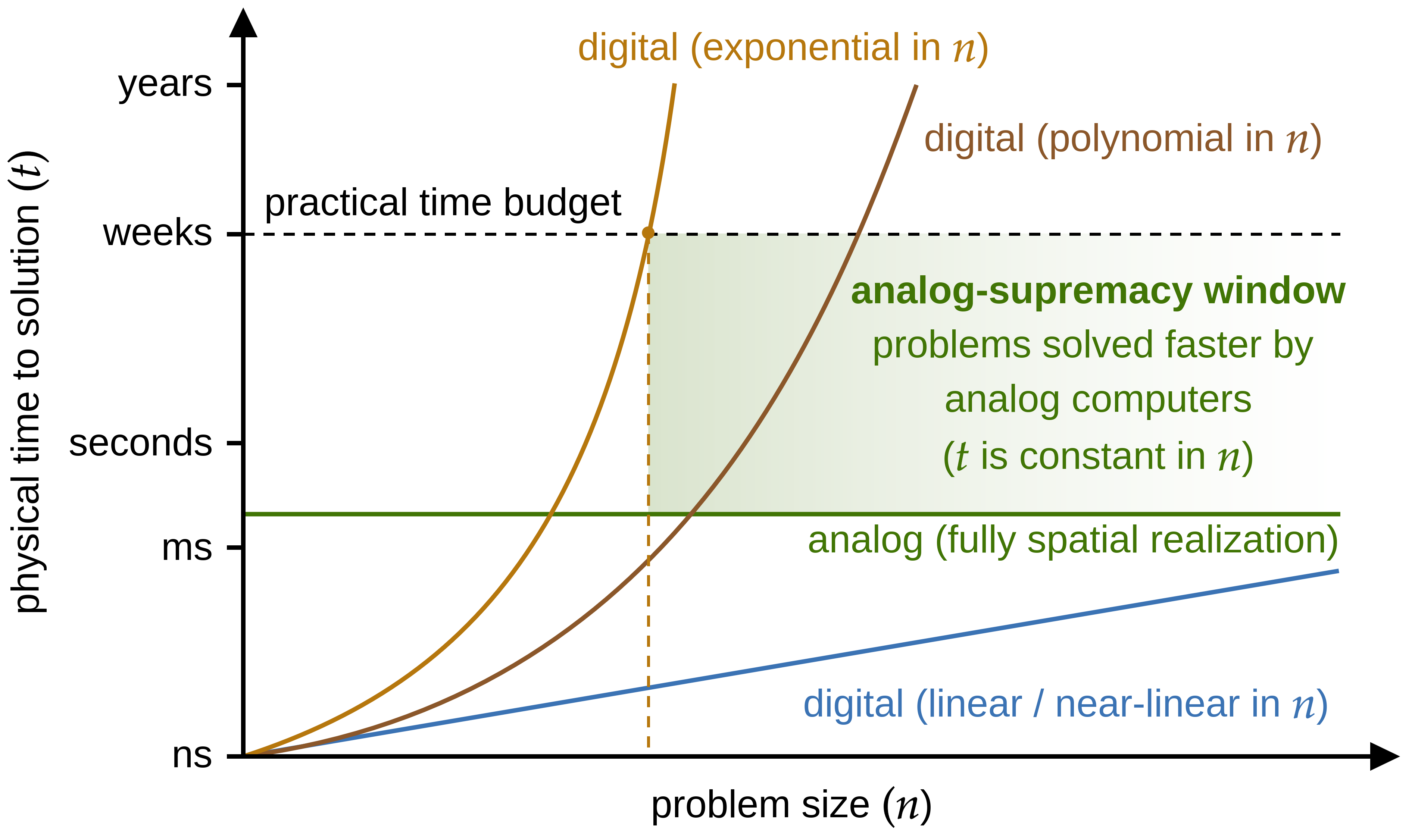}
    \caption{Conceptualization of the concept of ``analog supremacy''. The vertical axis is illustrative and indicates typical timescales for digital and analog computing.}
    \label{fig:analog_supremacy}
\end{figure}

We posit that any usage of the term \emph{analog supremacy}, visualized in Fig.~\ref{fig:analog_supremacy}, should distinguish between two claims. \emph{Soft analog supremacy} denotes an advantage concretely demonstrated on a finite and fully specified task. For a given problem of given size, at a stated error tolerance, an analog or hybrid computer outperforms the best available digital implementation in latency, throughput, energy per result, etc., after all overheads have been accounted for. The comparisons pertaining to the REDAC system described in Sec.~\ref{sec:redac_platform} and the results on analog attention in Sec.~\ref{sec:ai_applications} reflect claims of this kind. In contrast, and also based on Sec.~\ref{sec:new_quest}, we maintain that claims of \emph{hard analog supremacy} should be more strictly reserved for analog machines that can solve an entire problem family to $p$ bits of accuracy, while one or several of the accounted resources scale more favorably in the problem size $n$ and in $p$ than for the best known digital realization.\cite{JOC2007,JournalACM2017,Icalp2024} Before concluding, it is worth briefly mentioning that analog and quantum supremacy are conceptually related, but reflect profoundly different physical situations.\cite{Arute2019QuantumSupremacy} Quantum computation uses interference and entanglement that take place in a Hilbert space of dimension $2^n$ for $n$ qubits. A classical circuit representing a general $n$-qubit state would require $\Theta(2^n)$ state variables, and would therefore not inherit the advantages of quantum scaling.\cite{Feynman1982Simulating} However, even a hypothetically fault-tolerant quantum computer would not render analog computing irrelevant, because quantum computers can be expected to have their main impact on specific problem classes, most notably in chemistry and materials science,\cite{hoefler2023disentangling} and their practical reach is further constrained in the current noisy intermediate-scale quantum (NISQ) era.\cite{preskill2018nisq} By contrast, analog computers are far more general-purpose, being applicable as long as a problem can be expressed in terms of dynamical equations.
%
%%%%
%
\section{\label{sec:conclusion}Conclusion}

The need for energy-efficient, high-performance-computing solutions across a wide range of scientific and applied contexts has led to a resurgence of interest in analog and hybrid computing.

Recent progress in both theoretical foundations and hardware development has made it possible both to employ analog hardware in specialized scientific-computing workflows (such as matrix--vector multiplication) and to design preliminary versions of general-purpose analog computers (GPACs) and corresponding software stacks to make them accessible to academia and industry. 
Analog computing has clear applications in the simulation of deterministic and stochastic dynamical systems, and it also has considerable potential in areas like edge computing and the energy-efficient implementation of artificial-intelligence systems.

Modern analog computing has a very bright present and an even brighter future, and we highlight several of its open questions:
\begin{itemize}
\itemsep-0.25em
\item{Scaling and benchmarking: For which types of dynamical systems do analog computers provide substantial speedups and energy savings over commonly used numerical solvers in digital systems? How can one build fair benchmark suites (\ie, collections of standardized test problems) to compare their performance?}
\item{Error analysis: In computational mathematics, there are very well-developed theories of error propagation and error handling.\cite{higham2002accuracy}
It is important to develop similar theories for analog computing. The typical sources of imprecisions and errors are different in analog and digital computers, and the mechanisms and effects of error propagation also differ in the two computing paradigms.}
\item{Software support: What steps are necessary to seamlessly integrate current and future analog computer hardware into our predominantly algorithmic world? How can these systems be made broadly accessible and useful to academia and industry? Additionally, how (if at all) should software design 
influence decision-making in the development of the next generation of analog computers? 
The joint development of software and hardware, which is known as ``software--hardware codesign''\cite{mehdi2025codesign},
is a major driver in algorithm-aware hardware and software development in the field of artificial intelligence.}
\item{Artificial intelligence and large language models: Which components of large language models and related artificial-intelligence systems are most suitable for implementation in special-purpose analog hardware? It is also important to investigate the use of analog and hybrid computers, such as systems of interconnected biology-inspired spiking and bursting ``neurons'' and other coupled oscillators, as novel approaches for artificial-intelligence applications.
Furthermore, can analog computers be integrated into neural networks, thereby providing physical grounding for training and inference? Such advances have the potential to be very useful for robotics.\cite{stradmann2024closing,hines2025compact} 
}
\item{Neural ordinary differential equations: How can one use analog implementations of neural ODEs\cite{doi:10.1126/sciadv.adr7571,doi:10.1126/science.aee6277} for dynamics reconstruction, control, and their applications\cite{bottcher2022ai,bottcher2025control,bottcher2026control}?}
\item{Mechanical and physical learning: How can one leverage mechanical and physical learning~\cite{stern2020supervised,guzman2026unsupervised} in analog-computing hardware?}
\item{Combinatorial-optimization benchmarks: Which benchmarks are needed to demonstrate meaningful advantages of analog implementations over existing digital combinatorial-optimization solvers?}
\end{itemize}
\begin{acknowledgements}
    This perspective article is a synthesis of the presentations and discussions that took place at \href{https://analog-computing-days.org/}{Analog Computing Days 2026}. 
    We thank the Frankfurt School of Finance \& Management for supporting the organization of the event and acknowledge support from the Hessian Center for Artificial Intelligence (hessian.AI), Germany.
    
\end{acknowledgements}

%%%%
\section*{Data availability statement}

The data that support the findings of this study are available from the corresponding author upon reasonable request. 

\section*{Statement on AI use}

MAP did not use any AI for his contributions. LB used ChatGPT 5.5 to generate initial versions of the TikZ code for Figs.~\ref{pic_decay_example} and \ref{fig:linear_system}. LB subsequently edited both figures manually. LW used Claude Opus 4.8 to review structure, grammar, and spelling in Sec.~\ref{sec:general_purpose}.
%%%

%
%
%
\section*{Competing interests}
BU, LW, DT, CJ, SO, and PP are employees of anabrid GmbH, a company that develops analog computers. The other authors declare no competing interests.
\bibliography{analog-v03,olivier,michael,paul,software,bhaskar,refs_subsection_redac,analog_computing_theory_references}
% Produces the bibliography via BibTeX.

\end{document}